\documentclass[a4paper, 11pt]{article}
\title{Deformations and Fourier-Mukai transforms}
\date{}
\author{Yukinobu Toda}
\usepackage{amsmath,amssymb}
\usepackage{array}
\usepackage{amscd}

\newcommand{\aA}{\mathcal{A}}
\newcommand{\bB}{\mathcal{B}}
\newcommand{\cC}{\mathcal{C}}

\newcommand{\eE}{\mathcal{E}}
\newcommand{\fF}{\mathcal{F}}
\newcommand{\gG}{\mathcal{G}}

\newcommand{\mM}{\mathcal{M}}

\newcommand{\oO}{\mathcal{O}}
\newcommand{\pP}{\mathcal{P}}
\newcommand{\qQ}{\mathcal{Q}}

\newcommand{\uU}{\mathcal{U}}

\newcommand{\xX}{\mathcal{X}}
\newcommand{\yY}{\mathcal{Y}}

\newcommand{\lr}{\longrightarrow}

\newcommand{\Hom}{\mathop{\rm Hom}\nolimits}

\newcommand{\dotimes}{\stackrel{\textbf{L}}{\otimes}}
\newcommand{\dR}{\mathbf{R}}
\newcommand{\dL}{\mathbf{L}}

\newcommand{\Pic}{\mathop{\rm Pic}\nolimits}

\newcommand{\id}{\textrm{id}}

\newcommand{\Ext}{\mathop{\rm Ext}\nolimits}
\newcommand{\Spec}{\mathop{\rm Spec}\nolimits}

\newcommand{\Coh}{\mathop{\rm Coh}\nolimits}
\newcommand{\QCoh}{\mathop{\rm QCoh}\nolimits}
\newcommand{\Mod}{\mathop{\rm Mod}\nolimits}

\newcommand{\im}{\mathop{\rm im}\nolimits}

\newcommand{\cneq}{\mathrel{\raise.095ex\hbox{:}\mkern-4.2mu=}}
\newcommand{\eqcn}{\mathrel{=\mkern-4.5mu\raise.095ex\hbox{:}}}

\newcommand{\hocolim}{\mathop{\rm hocolim}\nolimits}

\newcommand{\kura}{\mathop{\rm Def}\nolimits}
\newcommand{\Cone}{\mathop{\rm Cone}\nolimits}
\newcommand{\perf}{\mathop{\rm perf}\nolimits}

\newcommand{\codim}{\mathop{\rm codim}\nolimits}

\newcommand{\lt}{\longrightarrow}

\newtheorem{thm}{Theorem}[section]
\newtheorem{prop}[thm]{Proposition}
\newtheorem{lem}[thm]{Lemma}
\newtheorem{defi}[thm]{Definition}
\newtheorem{rmk}[thm]{Remark}
\newtheorem{cor}[thm]{Corollary}

\newtheorem{prop-defi}[thm]{Proposition-Definition}

\begin{document}
\maketitle
\begin{abstract}
The aim of this paper is twofold. First we give an explicit
construction of the infinitesimal deformations of the category 
$\Coh (X)$ of coherent sheaves on a smooth projective variety 
$X$. Secondly, we show that any Fourier-Mukai transform 
$\Phi \colon D^b (X)\to D^b (Y)$ extends 
to an equivalence between the derived categories of the deformed 
Abelian categories.
\end{abstract}

\section{Introduction}

Recent developments on derived categories, coming from Homological mirror 
symmetry~\cite{Kon} or birational geometry~\cite{Ka1}, motivate
the necessity to establish a good deformation theory
of derived categories. 
The general deformation theory of Abelian categories was previously 
studied in~\cite{VB}, and the $A_{\infty}$-deformations of triangulated 
categories were studied in~\cite{Bara}.
However
 these analysis in these papers does not address the relationship 
 between deformations and Fourier-Mukai transforms. 
 So the following 
question arises:
 
 ``How do deformations interact with Fourier-Mukai transforms?"  
 
In this paper we concentrate on the first order deformations of $\Coh (X)$,
 and answer the above question in this case. Here 
$X$ is a smooth projective variety and $\Coh (X)$ is an Abelian category of coherent sheaves on $X$. By the philosophy of Kontsevich~\cite{Kon},
the Hochschild cohomology $HH^{\ast}(X)$ should parameterize 
 deformations of derived categories. 
 The degree 2-part should consist of deformations of $\Coh (X)$, since $HH^{2}(X)$ contains $H^1 (X,T_X)$ 
(deformations of complex structures) as a direct summand. The famous HKR-isomorphism says that $N$-th Hochschild cohomology is isomorphic to the 
direct sum $HT^N (X)\cneq \oplus _{p+q=N}H^p (X,\wedge ^q T_X)$. So there should be $\mathbb{C}[\varepsilon]/(\varepsilon ^2)$-linear Abelian category 
$\Coh (X,u)$ for $u\in HT^2 (X)$. Roughly the goals of this paper can
be summarized as follows.
\begin{itemize}
\item Give an explicit construction of $\mathbb{C}[\varepsilon]/(\varepsilon ^2)$-linear Abelian category $\Coh (X,u)$. 
\item Understand the behavior of the deformed triangulated category
$$D^b(X,u)\cneq D^b (\Coh (X,u))$$ under Fourier-Mukai transform $\Phi \colon D^b (X)\to D^b (Y)$. 
\end{itemize}

Note that any $u\in HT^2 (X)$ can be written as a sum $\alpha +\beta +\gamma$, 
with $\alpha \in H^2 (X,\oO _X)$, $\beta \in H^1 (X,T_X)$, and $\gamma \in H^0 (X,\wedge ^2 T_X)$. Then $\beta$ corresponds to a  
deformation of $X$ as a scheme, $\gamma$ is a non-commutative deformation. 
We will introduce ``twisted" sheaves using $\alpha$, 
and define $\Coh (X,u)$ as a combination of these components.

Next we make the second goal more precise. 
 Let $X$ and $Y$ be smooth projective varieties such that 
there exists an equivalence $\Phi \colon D^b (X)\to D^b (Y)$. Then 
we have an induced isomorphism of Hochschild cohomologies 
$\phi \colon HH^{\ast}(X) \to HH^{\ast}(Y)$. By combining 
$\phi$ with
HKR isomorphisms, we obtain the isomorphism $\phi _T \colon HT^2 (X)\to HT^2 (Y)$. 
Then the main theorem of this paper is the following:
\begin{thm}\label{main}
For $u\in HT^2 (X)$, let $v\cneq \phi _T (u)\in HT^2 (Y)$. Then there 
exists an equivalence 
$$\Phi ^{\dag}\colon D^b (X,u)\to D^b (Y,v),$$
such that the following diagram is 2-commutative.
$$\begin{CD}
D^b (X) @>{i_{\ast}}>> D^b (X,{u}) @>{\dL i^{\ast}}>> D^{-}(X)\\
@V{\Phi}VV @VV{\Phi ^{\dag}}V @VV{\Phi ^{-}}V\\
D^b (Y) @>{i_{\ast}}>> D^b (Y, v)@>{\dL i^{\ast}}>> D^{-}(Y).
\end{CD}$$
\end{thm}

By the above theorem, we can compare deformation theories under Fourier-Mukai transforms. One of the interesting point of Theorem~\ref{main} is that 
$\phi _T$ does not necessary 
preserve direct summands of $HT^2 (X)$. This indicates
$\Phi ^{\dag}$ may produce new interesting Fourier-Mukai dualities, 
for example dualities between usual commutative schemes and non-commutative 
schemes. Recently in the paper~\cite{OJT}, the equivalence $\Phi ^{\dag}$ 
of Theorem~\ref{main} has been 
 extended to infinite order deformations, when $X$ is an Abelian 
variety, $Y$ is its dual, and $\Phi$ is given by the Poincare line bundle. 
This result is giving a new kind of dualities via deformations, and 
it seems we will be able to find more examples of Fourier-Mukai equivalences
through deformation methods.  

\textbf{Acknowledgement}
 The author would like to express his profound gratitude to 
Professor Yujiro Kawamata, 
for many valuable comments, and warmful encouragement. 
The author also would like to thank T.Bridgeland for informing 
the author on the paper~\cite{Ne}.

\section{Hochschild cohomology and derived category}
Let $X$ be a smooth projective variety over $\mathbb{C}$ and
$\Delta _X \subset X\times X$ be 
a diagonal. We write $\Delta _X$ as $\Delta$ if it causes no confusion. 
In this section we recall the definitions of Fourier-Mukai transform, Hochschild cohomology and their properties. 
\begin{defi} Let $X$ and $Y$ be smooth projective varieties and take 
$\pP \in D^b (X\times Y)$. 
Let $p_i$ be projections from $X\times Y$ onto the corresponding factors. 
We define $\Phi _{X\to Y}^{\pP}$ 
as the following functor:
$$\Phi _{X\to Y}^{\pP}\cneq \dR p_{2\ast}(p_1 ^{\ast}(\ast)\dotimes \pP)\colon D^b (X)\to D^b (Y).$$
$\Phi _{X\to Y}^{\pP}$ is called integral transform with kernel $\pP$.
 If $\Phi _{X\to Y}^{\pP}$ gives an equivalence, then 
it is called a Fourier-Mukai transform. 
\end{defi}
The following theorem is fundamental in studying derived categories. 
\begin{thm}[Orlov~\cite{Or1}]
Let $\Phi \colon D^b (X) \to D^b (Y)$ be an exact functor. Assume that $\Phi$ is fully faithful and has a right adjoint.
Then there exists an object $\mathcal{P} \in D^b (X\times Y)$ such that $\Phi$ is isomorphic to the functor 
$\Phi_{X\to Y} ^{\mathcal{P}}$. Moreover $\pP$ is uniquely determined up to isomorphism. 
\label{O1}
\end{thm} 
Next we recall the Hochschild cohomology of the structure sheaf, given 
in~\cite{Kon}. 
\begin{defi}We define $HH^N(X)$ and $HT^N(X)$ as follows:
\begin{align*}HH^N(X) &\cneq \Hom _{X\times X}
(\oO _{\Delta}, \oO _{\Delta}[N]), \\
HT^N (X)&\cneq \bigoplus _{p +q=N}H^p (X, \bigwedge ^q T_X).\end{align*}
Here $\Hom$ is a morphism in $D^b (X\times X)$. $HH^{\ast}(X)$ is called Hochschild cohomology. 
\end{defi}
Note that the object $\fF \in D^b (X\times X)$ gives a functor 
$\Phi _{X\to X}^{\fF}$, and the morphism $\fF \to \gG$ gives a natural transformation $\Phi _{X\to X}^{\fF}\to \Phi _{X\to X}^{\gG}$. 
In this sense, Hochschild cohomology is a natural transformation $\id _X \to [N]$. 
But as in~\cite{Cal}, we can not consider $D^b (X\times X)$ as the category of functors precisely. 
(The map from the morphisms in $D^b (X\times X)$ to the natural transformations is not injective in general.)
However we can show the several properties of derived categories concerning $D^b (X\times X)$, for example categorical invariance of 
Hochschild cohomology, as if it is a category of functors. Since the natural transformations are categorical, Hochschild cohomology 
should be categorical invariant. In fact we have the following 
theorem in~\cite{Cal}. 

\begin{thm}[Caldararu~\cite{Cal}]\label{hh}
Let $X$ and $Y$ be smooth projective varieties such that there exists an equivalence $\Phi \colon D^b (X)\to D^b (Y)$. Then $\Phi$ induces an 
isomorphism $\phi \colon HH^{\ast}(X)\to HH^{\ast}(Y)$. 
\end{thm}
\textit{Outline of the proof}.
We will give the outline of the Caldararu's proof. 
Let $\pP \in D^b (X\times Y)$ be a kernel of $\Phi$, and
$\eE \in D^b (X\times Y)$ be a kernel of $\Phi ^{-1}$. Let 
$p_{ij}\colon X\times X\times Y\times Y \to X\times Y$ be projections onto corresponding factors.
Caldararu~\cite{Cal} showed that the functor with kernel= $p_{13}^{\ast}\pP \boxtimes p_{24}^{\ast}\eE \in D^b (X\times X\times Y\times Y)$,
$$\Phi _{X\times X \to Y\times Y}^{p_{13}^{\ast}\pP \boxtimes p_{24}^{\ast}\eE}\colon D^b (X\times X)\to D^b (Y\times Y)$$
gives an equivalence which takes $\oO _{\Delta _X}$ to $\oO _{\Delta _Y}$.
This equivalence implies the theorem immediately. $\quad$ q.e.d

\hspace{5mm}

Next we can compare $HH^{\ast}(X)$ and $HT^{\ast}(X)$. Hochschild cohomology is useful since its definition is categorical. But it is difficult
 to write down Hochschild cohomology classes
  explicitly. In calculating Hochschild cohomology, we decompose it
into direct sums of sheaf cohomologies of tangent bundles. The following 
theorem is due to Hochschild-Kostant-Rosenberg~\cite{HKR}
, Kontsevich~\cite{Kon}, Swan~\cite{Sw}, and Yekutieli~\cite{Ye}.
\begin{thm}\label{HKR}
There exists an isomorphism,
 $$I_{HKR} \colon HT^{\ast}(X)\to HH^{\ast}(X).$$
\end{thm}
\textit{Outline of the proof}. 
Note that $HH^N (X) \cong \Hom _X(\dL \Delta ^{\ast}\oO _{\Delta}, \oO _{X}[N])$ by adjunction. 
Let $\oO _X ^{\otimes i}\in \Mod (\oO _X)$ be the sheaf associated to the  
following presheaf:
$$U\subset X \longmapsto \Gamma (U,\oO _X)^{\otimes i}.$$
Here $\otimes$ is over $\mathbb{C}$, and $\oO _X$-module structure
on $\oO _X ^{\otimes i}$ is 
given by 
$$a\cdot (x_0 \otimes x_1 \otimes \cdots \otimes x_i)\cneq ax_0 \otimes x_1
\otimes \cdots \otimes x_i,$$
for $a,x_k \in \oO _X$. 
Let $d^i \colon \oO _X ^{\otimes(i+1)}\to \oO _X ^{\otimes i}$ be 
\begin{multline*}
d^{i}(x_0 \otimes \cdots \otimes x_i )=\sum _{k=0}^{i-1} (-1)^k x_0 \otimes \cdots \otimes x_k x_{k+1} \otimes \cdots \otimes x_i \\
+(-1)^i x_0 x_i \otimes x_1 \otimes \cdots 
\otimes x_{i-1}.\end{multline*}
Then we have the complex of $\oO _{X}$-modules:
$$\cC _X \cneq (\to \oO _X ^{\otimes (i+1)}\stackrel{d^i}{\to}\oO _X ^{\otimes i}\to \cdots \to \oO _X \to 0).$$
By \cite{Ye}, we have an explicit quasi-isomorphism 
$\cC _X \stackrel{\sim}{\longrightarrow}\dL \Delta ^{\ast}\oO _{\Delta}$
in $D(\Mod (\oO _X))$. 
Yekutieli~\cite{Ye} describes this 
isomorphism by building a resolution
using the formal neighborhood $X\subset X\times X\times \cdots \times X$. 
On the other hand, we have the following quasi-isomorphism $\cC _X
\to \oplus _{p\ge 0}\Omega _X ^p [p]$:
$$\begin{CD}
@>>> \oO _X ^{\otimes (i+1)}@>{d^{i}}>> \cdots @>{d^{2}}>>
\oO _X \otimes _{\mathbb{C}}\oO _X @>{0}>> \oO _{X} @>>> 0 \\ 
@. @V{I ^i}VV  @.  @V{I ^1}VV  @V{I ^0}VV \\
@>>> \Omega ^i _X @>{0}>> \cdots @>{0}>> \Omega _X @>{0}>> \oO _X @>>>0.
\end{CD}$$
Here $I ^i \colon \oO _X ^{\otimes (i+1)}\to \Omega _X ^i$ is given by 
$$I ^i (x_0 \otimes \cdots \otimes x_i)=x_0 \cdot dx_1 \wedge \cdots \wedge dx_i.$$
One can consult~\cite{Lo} for the detail.
 Consequently we get the quasi-isomorphism, 
$I \colon \dL \Delta ^{\ast}\oO _{\Delta}\stackrel{\sim}{\longrightarrow}
\oplus _{p\ge 0}\Omega _X ^p [p]$. 
Therefore we have the following isomorphism:
$$\Hom _X (\oplus _{p\ge0}\Omega _X ^p [p], \oO _X [N]) \stackrel{I}{\to}\Hom _X (\dL \Delta ^{\ast}\oO _{\Delta}, \oO _X [N]).$$ The left hand side is $HT^N (X)$ and the right hand side is $HH^N (X)$. $\quad$ q.e.d

                                                                                                       \vspace{5mm}
                                                                                   
$I_{HKR}$ is called HKR(Hochschild-Kostant-Rosenberg)-isomorphism. In the rest of this paper we write $I_{HKR}$ as $I_X$. 
Assume that $X$ and $Y$ are related by some Fourier-Mukai transform $\Phi \colon D^b (X)\to D^b (Y)$. 
By combining the isomorphisms $I_X$, $I_Y$ and $\phi$, we have the isomorphism:
$$\phi _{T}\cneq I_Y ^{-1}\circ \phi \circ I_X \colon HT^{\ast}(X)\stackrel{\sim}{\longrightarrow}HT^{\ast}(Y).$$
In the following 2-sections, we will construct deformations of $\Coh (X)$ for $u\in HT^2 (X)$.

\section{Non-commutative deformations of affine schemes}
Let $R$ be a Noetherian commutative ring and $X=\Spec R$. 
In this section we will consider a sheaf $\aA$ 
of (not necessary commutative) algebras on $X$.
Let $\uU _X$ be the category whose 
objects consist of Zariski open subset of $X$, and $\aA$ be a sheaf of algebra on $X$. Recall that a sheaf $\mM$ of left 
$\aA$-modules is quasi-coherent if for each $x\in X$, there exists an open neighborhood $U$ of $x$ and an 
exact sequence of left $\aA _U$-modules,
$$(\aA _U)^J \longrightarrow (\aA _U)^I \longrightarrow \mM _U \longrightarrow 0. $$
$\mM$ is coherent if the following conditions are satisfied:
\begin{itemize}
\item{$\mM$ is finitely generated, i.e. for every $x\in X$, there exists an open neighborhood $U$ of $x$ and a surjection 
$(\aA _U)^n \twoheadrightarrow M_U$.} 
\item{For every $U\in \uU _X$ and every $n\in \mathbb{Z}_{>0}$, and an arbitrary morphism of left $\aA _U$-modules 
$\phi \colon (\aA _U)^n \to \mM _U$, $\ker \phi$ is finitely generated. }
\end{itemize}
We denote  by $\Mod (\aA)$ the category of sheaves of left $\aA$-modules, by $\QCoh (\aA)$ full-subcategory of quasi-coherent sheaves, and  by $\Coh (\aA)$ coherent sheaves.
Of course it is well-known that if $\aA=\oO _X$, then quasi-coherent sheaf is written as $\widetilde{M}$ for some 
$R$-module $M$, and coherent sheaf is $\widetilde{M}$ for a 
finitely generated $R$-module $M$. 
We generalize these results to some non-commutative situations.
Let $\gamma$ be a bidifferential operator
$\gamma \colon R\times R \longrightarrow R.$
Using $\gamma$ we define a
 (not necessary commutative) ring structure on $R[\varepsilon]/(\varepsilon ^2)$ as follows:
$$(a+b\varepsilon)\ast _{\gamma}(c+d\varepsilon)\cneq ac+ (\gamma (a,c)+ad+bc)\varepsilon ,$$
and denote it by $R^{(\gamma)}$. 
Let $M$ be a left $R^{(\gamma)}$-module. Then the functor 
$$\uU _X  \ni U \longmapsto \oO _X (U) ^{(\gamma)}\otimes _{R^{(\gamma)}}M\in (\mbox{left } \oO _X (U) ^{(\gamma)} \mbox{-modules})$$
determines a presheaf of sets on $X$. Let $\widetilde{M}$ be the associated sheaf. 
We have a sheaf of rings $\oO _X ^{(\gamma)}\cneq \widetilde{R^{(\gamma)}}$ and $\widetilde{M}$ is a left 
$\oO _X ^{(\gamma)}$-module. 
Note that 
since $\oO _X (U) ^{(\gamma)}$ is right $R^{(\gamma)}$-left $\oO _X (U) ^{(\gamma)}$-module, 
$\oO _X (U) ^{(\gamma)}\otimes _{R^{(\gamma)}}M$ has a left 
$\oO _X (U) ^{(\gamma)}$-module structure.

As in the commutative case, we have the following lemma.
\begin{lem}\label{coh}
(1) For $f\in R$, $\widetilde{M}(U_f)=R_f ^{(\gamma)}\otimes _{R^{(\gamma)}}M$. In particular $\widetilde{M}(X)=M$ and 
$\oO _X ^{(\gamma)}(X)=R^{(\gamma)}$. \\
(2) $\widetilde{M}$ is a quasi-coherent $\oO _X ^{(\gamma)}$-module. \\
(3) The functor 
$$(\emph{left }R^{(\gamma)}\emph{-mod})\ni M \longmapsto \widetilde{M}\in \QCoh (\oO _X ^{(\gamma)})$$
gives an equivalence of categories.\\
(4) For $\fF \in \QCoh (\oO _X ^{(\gamma)})$, $\fF$ is coherent if and only if $M=\fF (X)$ is a finitely generated left $R^{(\gamma)}$-module.
\end{lem}
\textit{Proof}. Note that if we consider an $R$-module $N$ as left $R^{(\gamma)}$-module by the surjection 
$R^{(\gamma)}\twoheadrightarrow R$, then the action of $\oO _X ^{(\gamma)}$ on $\widetilde{N}$ descends to $\oO _X$, 
and $\widetilde{N}$ is a quasi-coherent $\oO _X$-module. \\
(1) It suffices to show $\widetilde{M}(X)=M$.
By the construction of $\widetilde{M}$, we have the natural morphism $M\to \widetilde{M}(X)$. Applying $\otimes _{R^{(\gamma)}}M$ to 
the surjection $R^{(\gamma)}\twoheadrightarrow R$, we obtain the exact sequence
$$0\longrightarrow \ker (r) \longrightarrow M \stackrel{r}{\longrightarrow}R\otimes _{R^{(\gamma)}}M \longrightarrow 0,$$
and the left action of $R^{(\gamma)}$ on $\ker (r)$ and $R\otimes _{R^{(\gamma)}}M$ descends to $R$. Therefore 
$\widetilde{\ker (r)}$ and $\widetilde{ R\otimes _{R^{(\gamma)}}M}$ are quasi-coherent $\oO _X$-modules. By applying $M\mapsto \widetilde{M}$ and taking 
global sections, we obtain the commutative diagram:
$$\begin{CD} 0 @>>> \ker (r) @>>> \widetilde{M}(X) @>>> R\otimes _{R^{(\gamma)}}M @>>> 0 \\
@.  @| @AAA @| @. \\
0 @>>> \ker (r) @>>> M @>>> R\otimes _{R^{(\gamma)}}M @>>> 0.
\end{CD}$$
It is easy to check that the multiplicative set 
$S=\{ f^{\ast _{\gamma}n} \} _{n\ge 0} \subset R^{(\gamma)}$ satisfies the right and left Ore localization conditions, and 
$R_f ^{(\gamma)}$ is a localization $S^{-1}R^{(\gamma)}$. Therefore the functor $M\mapsto \widetilde{M}$ is an exact functor. 
Moreover since $H^1 (X,\widetilde{\ker r})=0$, the top diagram is exact. By the 5-lemma, we have the isomorphism $M\to \widetilde{M}(X)$. \\
(2) Since $M\mapsto \widetilde{M}$ is an exact functor, 
we have an exact sequence 
$$\widetilde{R^{(\gamma)J}} \longrightarrow \widetilde{R^{(\gamma)I}} \longrightarrow \widetilde{M} \longrightarrow 0.$$
(3) Take $\fF \in \QCoh (\oO _X ^{(\gamma)})$. Applying $\otimes _{\oO _X ^{(\gamma)}}\fF$ to the exact sequence 
$$0\longrightarrow \varepsilon \oO _X \longrightarrow \oO _X ^{(\gamma)} \longrightarrow \oO _X \longrightarrow 0, $$
we can easily see that $\fF$ is given as an extension of quasi-coherent $\oO _X$-modules. Therefore the problem is reduced to the 
following lemma:

\begin{lem}\label{fullyfaith}
 Let $D^b (R^{(\gamma)})$ be the bounded derived category of left $R^{(\gamma)}$-modules, and
$\Mod (X,\gamma)\cneq \Mod (\oO _X ^{(\gamma)})$. 
The functor
$$  D^b (R^{(\gamma)})\ni M \longmapsto \widetilde{M}\in D^b (\Mod (X,\gamma))$$
is fully faithful.
\end{lem}
\textit{Proof}. 
Take $M,N \in D^b (R^{(\gamma)})\subset D(R^{(\gamma)})$ and we will show that 
$$\Hom _{D(R^{(\gamma)})}(M,N)\longrightarrow \Hom _{D(\Mod (X,\gamma))}(\widetilde{M}, \widetilde{N})$$
is an isomorphism. 
By taking a free resolution, we may assume $M$ is a bounded above complex of free $R^{(\gamma)}$-modules. Let 
$M_k \cneq \sigma _{\ge -k}M$. Here $\sigma _{\ge -k}$ denotes the
 stupid truncation. 
Now we have a sequence of complexes $\to M_k \to M_{k+1} \to \cdots$ and if we take the homotopy colimit (cf.~\cite{Ne})
$$\oplus _k M_k \stackrel{s-\id}{\longrightarrow} \oplus _k M_k \longrightarrow \hocolim (M_k) \longrightarrow \oplus _k M_k [1],$$
then there exists a quasi-isomorphism $\hocolim (M_k)\to M$. Here $s$ is the shift map, whose coordinates are the natural maps $M_k \to M_{k+1}$.
Therefore we may assume $M$ is a finite complex of free $R^{(\gamma)}$-modules. Again by taking stupid truncations, we may assume 
$M=R^{(\gamma)}$. Since $N$ is bounded, we may assume $N=N'[k]$ for some left $R^{(\gamma)}$-module $N'$. Now it suffices to show that the 
map 
$$\Hom _{D(R^{(\gamma)})}(R^{(\gamma)}, N'[k])\longrightarrow \Hom _{D(\Mod (X,\gamma))}(\oO _X ^{(\gamma)}, \widetilde{N}'[k])$$
is an isomorphism. If $k<0$, then both sides are zero. If $k=0$, then both sides are $N'$. If $k>0$, then the left hand side is zero, so it suffices to 
show $H^k (X, \widetilde{N}')=0$ for $k>0$. But since $\widetilde{N}'$ is an extension of quasi-coherent $\oO _X$-modules, 
$H^k (X, \widetilde{N}')=0$ for $k>0$. $\quad$ q.e.d

\hspace{5mm}

(4) First we
 check that a submodule of a finitely generated $R^{(\gamma)}$-module
 is also finitely generated. 
In fact let $M$ be a finitely generated $R^{(\gamma)}$-module, and $N\subset M$ be a submodule. Then 
we have 
the natural morphism $g\colon N \to R\otimes _{R^{(\gamma)}}M$. 
It is enough to check that $\ker (g)$ and $\im (g)$ are finitely generated 
$R^{(\gamma)}$-modules. 
Note that we have $\ker (g)\subset \varepsilon M$ and 
$\im (g) \subset R\otimes _{R^{(\gamma)}}M$. Since $R$ is Noetherian
and $\varepsilon M$, $R\otimes _{R^{(\gamma)}}M$ are both 
finitely generated $R$-modules, it follows that $\ker (g)$ and $\im (g)$
are both finitely generated $R$-modules. Thus in particular these are 
finitely generated $R^{(\gamma)}$-modules via the surjection 
$R^{(\gamma)}\twoheadrightarrow R$.  

Using this fact, we can 
see $\widetilde{M}$ for a finitely generated left $R^{(\gamma)}$-module $M$ 
is coherent. On the other hand, take $\fF \in \Coh (\oO _X ^{(\gamma)})$. Then by (3),
$\fF$ can be written as $\fF =\widetilde{M}$ for some left 
$R^{(\gamma)}$-module $M$. Since $\fF$ is given by an extension of coherent $\oO _X$-modules, 
$M$ is a finitely generated left $R^{(\gamma)}$-module. $\quad$ q.e.d

\hspace{5mm}

For a full subcategory $\cC \subset \Mod (X,\gamma)$, 
let $D^b _{\cC}(\Mod (X,\gamma))$ denote the full subcategory of $D^b (\Mod (X,\gamma))$ whose objects have cohomologies contained in $\cC$. 
As a corollary, we obtain the following:
\begin{cor}\label{ne} There exist equivalences, 
$$D^b (R^{(\gamma)})\stackrel{\sim}{\longrightarrow}D^b _{\QCoh}(\Mod (X,\gamma)), \quad
D^b _f(R^{(\gamma)})\stackrel{\sim}{\longrightarrow}D^b _{\Coh}(\Mod (X,\gamma)).$$
Here $D^b _f(R^{(\gamma)})$ is a derived category of finitely generated left $R^{(\gamma)}$-modules.
\end{cor}
\textit{Proof}. We have proved the full faithfulness 
in Lemma~\ref{fullyfaith}. Since an object of $\QCoh (\oO _X ^{(\gamma)})$ is written as $\widetilde{M}$ for a
left $R^{(\gamma)}$-module $M$, the image from the left hand side generates 
the right hand side. $\quad$ q.e.d

\begin{rmk} In general we can show the unbounded case of the above corollary as in~\cite{Ne}.
 Here we gave a proof of bounded case for the sake of simplicity. 
For the details, the reader should refer to~\cite{Ne}. 
\end{rmk}

\section{Infinitesimal deformations of $\Coh (X)$}
From this section on, we will assume that
 $X$ is a smooth projective variety over $\mathbb{C}$. 
The aim of this section is to construct the 
first order deformations of $\Coh (X)$.
First we begin with the general situation.
Let us take an affine open cover $X=\cup _{i=1}^N U_i$, and denote by $\mathfrak{U}$ this open 
cover. 
Let $U_{i_0 \cdots i_p}\cneq U_{i_0}\cap \cdots \cap U_{i_p}$, and $j_{i_0 \cdots i_p}\colon U_{i_0 \cdots i_p}\hookrightarrow X$ be
open immersions. 
For a sheaf $\fF$ on $X$, let $\cC ^p (\mathfrak{U},\fF)$, $\mathbf{C}^p (\mathfrak{U},\fF)$
be
$$\cC ^p (\mathfrak{U},\fF)\cneq \prod _{i_0 \cdots i_p}j_{i_0 \cdots i_p \ast}j_{i_0 \cdots i_p}^{\ast}\fF, \quad 
\mathbf{C}^p (\mathfrak{U}, \fF)\cneq \prod _{i_0 \cdots i_p}\Gamma (U_{i_0 \cdots i_p}, \fF).$$
Let us consider a sheaf of algebras $\aA$ on $X$ and its center $Z(\aA)$. Take $\tau \in H^2 (X,Z(\aA)^{\times})$. Then 
$\tau$ is represented by a $\check{\mbox{C}}$ech cocycle 
$\tau=\{ \tau _{i_0 i_1 i_2}\} \in \mathbf{C}^2 (\mathfrak{U}, Z(\aA)^{\times}).$ We define the category $\Mod (\aA,\tau)$ as follows:

\begin{defi} We define $\Mod (\aA, \tau)$ as an Abelian category of $\tau$-twisted left $\aA$-modules. Namely objects of 
$\Mod (\aA, \tau)$ are collections
$$\fF=(\{ \fF _i \}_{1\le i \le N}, \phi _{i_0 i_1}),$$
where $\fF _i \in \Mod (\aA |_{U_i})$ and $\phi _{i_0 i_1}$ 
are isomorphisms
$$\phi _{i_0 i_1}\colon \fF _{i_0}|_{U_{i_0 i_1}}\stackrel{\cong}{\longrightarrow}\fF _{i_1}|_{U_{i_0 i_1}}$$
as left $\aA |_{U_i}$-modules. These data must satisfy the equality
$$\phi _{i_2 i_0}\circ \phi _{i_1 i_2}\circ \phi _{i_0 i_1} =\tau _{i_0 i_1 i_2} \cdot \emph{id} _{\fF _0}.$$
\end{defi}
We say $\fF \in \Mod (\aA, \tau)$ is quasi-coherent if $\fF _i \in \QCoh (\aA |_{U_i})$, and coherent if $\fF _i \in \Coh (\aA |_{U_i})$. 
We denote by $\QCoh (\aA, \tau)$ the category of quasi-coherent $\tau$-twisted left $\aA$-modules, and by $\Coh (\aA, \tau)$ coherent twisted 
sheaves. 
\begin{lem}
Up to equivalence, the categories
$\Mod (\aA, \tau)$, $\QCoh (\aA, \tau)$, $\Coh (\aA, \tau)$ are independent of choices of $\mathfrak{U}$ 
and $\check{\mbox{C}}$ech representative of $\alpha$. \end{lem}
\textit{Proof}.
The proof is easy and left it to the reader. $\quad$ q.e.d

\subsection*{Fundamental properties and operations on $\Mod (\aA, \tau)$}
\begin{itemize}
\item  \textit{$j^{\ast}$, $j_{\ast}$, $j_{!}$ for an open immersion $j\colon U\hookrightarrow X$} \end{itemize}
Let $j\colon U\hookrightarrow X$ be an open immersion. We have the obvious functors:
\begin{align*}
j^{\ast}\colon \Mod (\aA, \tau) &\longrightarrow \Mod (\aA |_{U}, \tau |_{U}),\\ j_{\ast}, j_{!}\colon 
\Mod (\aA |_{U}, \tau |_{U})&\longrightarrow \Mod (\aA, \tau).\end{align*}
$j_{\ast}$ is right adjoint of $j^{\ast}$, and $j_{!}$ is left adjoint of $j^{\ast}$.
For $\fF =(\{\fF _{i_0} \} , \phi _{i_0 i_1})\in \Mod (\aA |_{U}, \tau _{U})$, 
with $\fF _{i_0}\in \Mod (\aA|_{U\cap U_{i_0}})$, $j_{\ast}\fF$ and $j_{!}\fF$ are given by 
$$j_{\ast}(\fF)_{i_0} \cneq (j|_{U\cap U_{i_0}})_{\ast}\fF _{i_0}, \quad j_{!}(\fF )_{i_0}\cneq (j|_{U\cap U_{i_0}})_{!}\fF _{i_0}.$$
Here 
$$(j|_{U\cap U_{i_0}})_{!}\colon \Mod (\aA|_{U\cap U_{i_0}})\longrightarrow
\Mod (\aA|_{U_{i_0}})$$
is extension by zero.

\begin{itemize}
\item  \textit{Tensor product}\end{itemize} 
Let us take $\fF \in \Mod (\aA ^{op},\tau)$. Assume that the right action of the subalgebra $\bB \subset \aA$ on $\fF$ is centralized. Then we
have the functor,
$$\fF \otimes \ast \colon \Mod (\aA, \tau ')\longrightarrow \Mod (\bB, \tau \cdot \tau ').$$
In particular if $\bB$ is contained in the center of $\aA$, then we have the functor,
$$\otimes \colon \Mod (\aA ^{op}, \tau)\times \Mod (\aA, \tau ')\longrightarrow \Mod (\bB, \tau \cdot \tau ').$$

\begin{itemize}
\item \textit{Pull-back} \end{itemize}
Let $f\colon Y\to X$ be a morphism of varieties, and $\aA$, $\bB$ be sheaves of algebra on $X$ and $Y$. If 
there exists a morphism of algebras $f^{-1}\aA \to \bB$ which preserves their centers, then we have the pullback
$$f^{\ast}\colon \Mod (\aA, \tau) \longrightarrow \Mod (\bB, f^{\ast}\tau),$$
which takes $(\{ \fF _i \}, \phi _{i_0 i_1})$ to 
$(\{ \bB \otimes _{\aA}f^{-1}\fF _i \}, 1\otimes \phi _{i_0 i_1})$.

\begin{itemize}
\item \textit{Push-forward} \end{itemize}
In the same situation as above, we have a morphism of algebras $\aA \to f_{\ast}\bB$ which preserves their centers. 
We have the push-forward:
$$f_{\ast}\colon \Mod (\bB, f^{\ast}\tau)\longrightarrow \Mod (\aA, \tau).$$
Clearly $f_{\ast}$ is a right adjoint of $f^{\ast}$.

\begin{itemize}
\item \textit{Enough injectives and flats} \end{itemize}
\begin{lem}\label{en}(i) $\Mod (\aA,\tau)$ has enough injectives. \\
(ii) For every $A\in \Mod (\aA, \tau)$, there exists a flat object $P\in \Mod (\aA, \tau)$ and a surjection 
$P\twoheadrightarrow A$. Here we say $\fF =(\{\fF_i \}, \phi _{i_0 i_1})$ is flat if each $\fF _i$ is a flat $\aA _{U_i}$-module. 
\end{lem}
\textit{Proof}. (i) Take $A\in \Mod (\aA, \tau)$. Since $\Mod (\aA |_{U_i})$ has enough injective, there 
exists an injection $j^{\ast}A \hookrightarrow I_i$ for an injective object $I_i \in \Mod (\aA |_{U_i})$. 
Let $\widetilde{I}_i\cneq j_{\ast}I_i$. Then the composition
$$A\longrightarrow j_{\ast}j^{\ast}A\longrightarrow \prod _i \widetilde{I}_i $$
is an injection. Since $j_{\ast}$ is a right adjoint of $j^{\ast}$, $\prod _i \widetilde{I}_i$ is an injective object of 
$\Mod (\aA,\tau)$. \\
(ii) Take $A\in \Mod (\aA, \tau)$. We can take a surjection $P_i \twoheadrightarrow j^{\ast}A$ for flat $\oO _{U_i}^{(\gamma)}$-module
$P_i$. Let $\bar{P_i}\cneq j_{!}P_i$. Then the composition
$$\bigoplus _i \bar{P_i} \longrightarrow j_{!}j^{\ast}A \longrightarrow A$$
is surjective and $\bigoplus _i \bar{P_i}$ is flat. $\quad$ q.e.d

\hspace{5mm}

Let us take an element
$$u=(\alpha, \beta, \gamma)\in HT^2 (X)=H^2 (\oO _X)\oplus H^1 (T_X)\oplus H^0 (\wedge ^2 T_X).$$
First
we construct a sheaf $\oO _X ^{(\beta, \gamma)}$ of
 $\mathbb{C}[\varepsilon]/(\varepsilon ^2)$-algebras on $X$.
Note that we
 can consider $\gamma$ as a bidifferential operator $\oO _X \times \oO _X \to \oO _X$, and $\beta _{i_0 i_1}$ as a 
differential operator $\oO _{U_{i_0 i_1}}\to \oO _{U_{i_0 i_1}}$. 
As a sheaf $\oO _X ^{(\beta, \gamma)}$ is $\oO _X ^{(\beta _{i_0 i_1})}$, the kernel of the following morphism:
$$\oO _X \oplus \cC ^0 (\mathfrak{U}, \oO _X)\ni (a, \{b_i \})\stackrel{\delta ^{(\beta _{i_0 i_1})}}
{\longmapsto} -\beta _{i_0 i_1}(a)+\delta \{b_i \} \in \cC ^1 (\mathfrak{U}, \oO _X).$$ 
We define the product on $\oO _X \oplus \cC ^0 (\mathfrak{U}, \oO _X)$ by the formula: 
$$(a, \{b_i \})\ast _{\gamma}(c,\{d_i \})\cneq (ac, \{ ad_i +cb_i +\gamma (a,c)\} _i).$$
Then it is easy to see that $\oO _X ^{(\beta _{i_0 i_1})}$ is a subalgebra
of $\oO _X \oplus \cC ^0 (\mathfrak{U}, \oO _X)$, and denote by $\oO _X ^{(\beta, \gamma)}$ this sheaf of algebras. It is also easy to
check $\oO _X ^{(\beta, \gamma)}$ doesn't depend on the choices of $\mathfrak{U}$ and 
$\check{\mbox{C}}$ech representative of $\beta$.
Note that $\oO _X ^{(\beta, \gamma)}|_{U_i}\cong \oO _{U_i}^{(\gamma)}$ as a sheaf of algebra. 
Since $(1-\alpha _{i_0 i_1 i_2}\varepsilon)$ is contained in the center of $\oO _{U_{i_0 i_1 i_2}}^{(\gamma)}$, we have an 
element
$$\widetilde{\alpha}\cneq \{(1-\alpha _{i_0 i_1 i_2}\varepsilon)\} _{i_0 i_1 i_2} \in \mathbf{C}^2 (X, Z(\oO _X ^{(\beta, \gamma)})),$$
which is a cocycle. 
Let $\Mod (X,u)\cneq \Mod (\oO _X ^{(\beta, \gamma)}, \widetilde{\alpha})$, and define $\QCoh (X,u)$ and $\Coh (X,u)$ as above.

\hspace{5mm}

Now we can define $D^{\ast}(X,u)$ for $\ast =b, \pm, \emptyset$ as follows. 
\begin{defi}We define $\mathbb{C}[\varepsilon]/(\varepsilon ^2)$-linear triangulated category $D^{\ast}(X, u)$ as 
$$D^{\ast}(X,u)\cneq D^{\ast} (\Coh (X, u)),\qquad (\ast =b, \pm, \emptyset).$$
\end{defi}

As in~\cite{Ne}, we have the following proposition:
\begin{prop}\label{ma} There exist natural equivalences:
\begin{align*}
D^{\ast} (\QCoh (X,u))&\stackrel{\sim}{\longrightarrow} D_{\QCoh}^{\ast} (\Mod (X,u)), \\
D^{\ast}(X,u)&\stackrel{\sim}{\longrightarrow} D_{\Coh}^{\ast} 
(\Mod (X,u)), \end{align*}
for $\ast=b, \pm, \emptyset$.
\end{prop}
\textit{Proof}. The proof is the same as in~\cite{Ne}. Take an affine open cover $X=\cup _{i=1}^N U_i$. 
We use the induction on $N$ to prove the proposition, and the case of $N=1$ and $\ast =b$ 
has been proved in the previous section. 
$\quad$ q.e.d

\hspace{5mm}

Now we can construct transformations between derived categories.  Take two smooth projective varieties $X$ and $Y$, and 
$u=(\alpha, \beta, \gamma)\in HT^2(X)$, $v=(\alpha ',\beta ' ,\gamma ')\in HT^2 (Y)$. 
 For a perfect object (i.e. locally quasi-isomorphic to bounded complexes of free modules) 
$\pP ^{\dag}\in D^b _{\perf}(X\times Y, -p_1^{\ast}\check{u}+p_2 ^{\ast}v)$, we will construct a functor, 
$$\Phi ^{\dag}\colon D^b (X, {u})\longrightarrow D^b (Y,v).$$
Here $\check{u}\cneq (\alpha, -\beta, \gamma)$. Firstly take $\fF \in D^b(X,{u})$. Since we have a morphism of algebras
$$p_1 ^{-1}\oO _X ^{(\beta, \gamma)}\longrightarrow \oO _{X\times Y}^{(p_1 ^{\ast}\beta +p_2 ^{\ast}\beta ' , p_1 ^{\ast}\gamma -p_2 ^{\ast}\gamma ')},$$
we obtain the object 
\begin{align*}
p_{1}^{\ast}\fF &\in D^b(\Coh (\oO _{X\times Y}^{(p_1 ^{\ast}\beta +p_2 ^{\ast}\beta ' , p_1 ^{\ast}\gamma -p_2 ^{\ast}\gamma ')}, p_1 ^{\ast}\widetilde{\alpha})) \\
&\simeq D^b(\Coh (\oO _{X\times Y}^{(p_1 ^{\ast}\beta +p_2 ^{\ast}\beta ' , -p_1 ^{\ast}\gamma +p_2 ^{\ast}\gamma '),op}, p_1 ^{\ast}\widetilde{\alpha})).
\end{align*}
Now by Lemma~\ref{en}, we can define $\dotimes \pP ^{\dag}$. Since the right action of $p_2 ^{-1}\oO _Y ^{(\beta ', \gamma ')}$ 
on each term of $p_1 ^{\ast}\fF$ is centralized, we obtain the object, 
$$p_1 ^{\ast}\fF \dotimes \pP ^{\dag}\in D^b (\Mod (p_2 ^{-1}\oO _Y ^{(\beta ', \gamma ')}, p_2 ^{\ast}\widetilde{\alpha} ')).$$
(Since $\pP ^{\dag}$ is perfect, $\dotimes \pP ^{\dag}$ preserves boundedness.) Applying $\dR p_{2\ast}$, we obtain the object, 
$$\dR p_{2\ast}(p_1 ^{\ast}\fF \dotimes \pP ^{\dag})\in D^b (\Mod (\oO _Y ^{(\beta ', \gamma ')}, \widetilde{\alpha} ')).$$
If all the cohomologies $R^i p_{2\ast}(p_1 ^{\ast}\fF \dotimes \pP ^{\dag})$ are coherent, we can define $\Phi ^{\dag}$ as  
$$\Phi ^{\dag}(\fF)\cneq \dR p_{2\ast}(p_1 ^{\ast}\fF \dotimes \pP ^{\dag})\in D_{\Coh}^b (\Mod (\oO _Y ^{(\beta ' ,\gamma ')}, \alpha '))\simeq
D^b (Y,v),$$
by Lemma~\ref{ma}. In fact we have the following: 
\begin{lem}
For each $i$, the object $R^i p_{2\ast}(p_1 ^{\ast}\fF \dotimes \pP ^{\dag})$
 is coherent.
\end{lem}
\textit{Proof}.
Since there exists a distinguished triangle,
$$p_1 ^{\ast}\fF \dotimes \oO _{X\times Y} \dotimes \pP ^{\dag} \longrightarrow p_1 ^{\ast}\fF \dotimes \pP ^{\dag} \longrightarrow p_1 ^{\ast}\fF 
\dotimes \oO _{X\times Y} \dotimes \pP ^{\dag}$$
in $D^b (\Mod (p_2 ^{-1}\oO _Y ^{(\beta ', \gamma ')}, p_2 ^{\ast}\alpha '))$,
it suffices to show that
$R^i p_{2\ast}(p_1 ^{\ast}\fF \dotimes \oO _{X\times Y} \dotimes \pP ^{\dag})$
is coherent. But 
since $H^q (p_1 ^{\ast}\fF \dotimes \oO _{X\times Y} \dotimes \pP ^{\dag})$ are coherent $\oO _{X\times Y}$-modules, the existence of 
a first quadrant spectral sequence
$$E_2 ^{p,q}\cneq R^p p_{2\ast}H^q (p_1 ^{\ast}\fF \dotimes \oO _{X\times Y} \dotimes \pP ^{\dag}) 
\Rightarrow R^{p+q}p_{2\ast}(p_1 ^{\ast}\fF \dotimes \oO _{X\times Y} \dotimes \pP ^{\dag})$$
shows $R^i p_{2\ast}(p_1 ^{\ast}\fF \dotimes \oO _{X\times Y} \dotimes \pP ^{\dag})$ is coherent. $\quad$ q.e.d

\hspace{5mm}

Since we have a morphism of algebras $i\colon \oO _{X}^{(\beta,\gamma)}\to \oO _X$, we have functors:
$$i_{\ast}\colon \Coh (X) \to \Coh (X,u), \qquad i^{\ast}\colon \Coh (X,u)\to \Coh (X).$$
Passing to derived categories and using Proposition~\ref{ma}, we obtain the derived functors:
$$i_{\ast}\colon D^b (X)\to D^b (X,u), \qquad \dL i^{\ast} \colon D^b (X,u)\to D^{-}(X).$$
Note that an equivalence $\Phi \colon D^b(X)\to D^b(Y)$ extends to an equivalence $\Phi ^{-}\colon D^{-}(X)\to D^{-}(Y)$, using the same kernel with $\Phi$. 
Now we can state our main theorem. 
\begin{thm} \label{mt}
Let $X$ and $Y$ be smooth projective varieties such that there exists an equivalence of derived categories $\Phi \colon D^b (X)\to D^b (Y)$. 
Take $u\in HT^2 (X)$ and $v\cneq \phi _T (u)\in HT^2 (Y)$.
Then there exists an object $\pP ^{\dag}\in D^b _{\perf}(X\times Y, -p_1 ^{\ast}\check {u}+p_2 ^{\ast}v)$ such that the associated functor
$$\Phi ^{\dag}\colon D^b (X,{u})\longrightarrow D^b (Y, v), $$
gives an equivalence. Moreover the following diagram is 2-commutative. 
$$\begin{CD}
D^b (X) @>{i_{\ast}}>> D^b (X,{u}) @>{\dL i^{\ast}}>> D^{-}(X)\\
@V{\Phi}VV @VV{\Phi ^{\dag}}V @VV{\Phi ^{-}}V\\
D^b (Y) @>{i_{\ast}}>> D^b (Y, v)@>{\dL i^{\ast}}>> D^{-}(Y).
\end{CD}$$
\end{thm}

\section{Atiyah classes and FM-transforms} 
In this section we will analyze Atiyah classes of kernels of Fourier-Mukai transforms, and give the preparation 
for the proof of the main theorem. 
Firstly let us recall the universal Atiyah class. 
Let $X$ be a smooth projective variety and $\Delta$ be a diagonal or diagonal embedding. We write $\Delta$ as $\Delta _X$ when needed. 
Let $I_{\Delta}\subset \oO _{X\times X}$ be an ideal sheaf of $\Delta$. 
Consider the exact sequence,
$$0 \longrightarrow I_{\Delta}/I_{\Delta}^2 \longrightarrow \oO _{X\times X}/I_{\Delta}^2 \longrightarrow \oO _{\Delta}\longrightarrow 0.\quad \cdots (\star)$$
\begin{defi}The universal Atiyah class 
$$a_{X}\colon \oO _{\Delta}\longrightarrow \Delta _{\ast}\Omega _X [1]$$
is the extension class of the exact sequence $(\star)$. 
\end{defi}
Consider the composition
$$\oO _{\Delta}\stackrel{a_{X}}{\longrightarrow}\Delta _{\ast}\Omega _X [1]\stackrel{a_{X}\otimes p_2 ^{\ast}\Omega _X}{\longrightarrow}
\Delta _{\ast}\Omega _X ^{\otimes 2}[2]\longrightarrow \cdots
\longrightarrow \Delta _{\ast}\Omega _X ^{\otimes i}[i].$$
By composing anti-symmetrization $\epsilon \colon \Omega _X ^{\otimes i}\to \Omega _X ^i$, we get a morphism
$$a_{X,i}\colon \oO _{\Delta}\longrightarrow \Delta _{\ast}\Omega _X ^i [i].$$
\begin{defi} The exponential universal Atiyah class is a morphism
$$\exp (a)_X \cneq \bigoplus _{i\ge 0}a_{X,i}\colon \oO _{\Delta}\longrightarrow \bigoplus _{i\ge 0}\Delta _{\ast}\Omega _X ^i [i].$$
Here $a_{X,0}=\emph{id}$. 
\end{defi}
Caldararu~\cite{Cal2} showed the following:
\begin{prop}[Caldararu~\cite{Cal2}]\label{ati}
$\exp (a)_X$ is equal to the composition
$$\oO _{\Delta}\longrightarrow \Delta _{\ast}\dL \Delta ^{\ast}\oO _{\Delta}\stackrel{\Delta _{\ast}I}{\longrightarrow}\bigoplus _{i\ge 0}\Delta _{\ast}\Omega _X ^i [i].$$
Here $\oO _{\Delta}\to \Delta _{\ast}\dL \Delta ^{\ast}\oO _{\Delta}$ is an adjunction, and $I$ is a morphism which appeared in the proof of Theorem~\ref{HKR} 
\end{prop}
By the above proposition, HKR-isomorphism is nothing but the following morphism
$$HT^{\ast}(X)\ni u \longmapsto \Delta _{\ast}u \circ \exp (a)_X \in HH^{\ast}(X).$$

Next, let us recall the Atiyah class and exponential Atiyah class for an object $\pP \in D^b (X)$. 
By applying $\dR p_{2\ast}(p_1 ^{\ast}\pP \dotimes \ast)$ to the exact sequence $(\star)$, we obtain the distinguished triangle, 
$$\pP \otimes \Omega _X \longrightarrow \dR p_{2\ast}\left( p_1 ^{\ast}\pP \dotimes \oO _{X\times X}/I_{\Delta}^2 \right)
\longrightarrow \pP \longrightarrow \pP \otimes \Omega _X [1]. \quad \cdots (\star _{\pP})$$
\begin{defi}
The Atiyah class $a(\pP )\in \Ext _X ^1 (\pP , \pP \otimes \Omega _X)$ is a morphism
$$a(\pP)\colon \pP \longrightarrow \pP \otimes \Omega _X [1]$$
in the distinguished triangle $(\star _{\pP})$.
\end{defi}
As in the exponential universal Atiyah class, let us take the composition,
$$a(\pP)\circ \cdots a(\pP)\colon \pP \longrightarrow \pP \otimes \Omega _X [1]\longrightarrow  \cdots \longrightarrow \pP \otimes \Omega _X ^{\otimes i} [i].$$By composing $\epsilon \colon \Omega _X ^{\otimes i}\to \Omega _X ^i$, we get the morphism,
$$a(\pP )_i \colon \pP \longrightarrow \pP \otimes \Omega _X ^i [i].$$
\begin{defi}The exponential Atiyah class of $\pP$ is a morphism
$$\exp a(\pP)\cneq \bigoplus _{i\ge 0}a(\pP)_i \colon \pP \longrightarrow \bigoplus _{i\ge 0}\pP \otimes \Omega _X ^i [i].$$
Here $a(\pP )_0=\emph{id}$.
\end{defi}
Now let us consider two smooth projective varieties $X$ and $Y$, and an equivalence of 
derived categories $\Phi \colon D^b (X)\to D^b (Y)$. Let $\pP \in D^b (X\times Y)$ be a kernel of $\Phi$. By Theorem~\ref{hh}, $\Phi$ induces the isomorphism
$\phi \colon HH^{\ast}(X)\to HH^{\ast}(Y)$. 
We have the following proposition:
\begin{prop}\label{co}
$\phi$ factors into the isomorphisms:
$$HH^{\ast}(X)\stackrel{\sim}{\longrightarrow} \Ext _{X\times Y}^{\ast}(\pP ,\pP) \stackrel{\sim}{\longrightarrow}HH^{\ast}(Y).$$

\end{prop}
\textit{Proof}. 
Let $p_{ij}$ be projections from $X\times Y\times Z$ onto corresponding factors. 
For $a\in D^b (X\times Y)$ and $b\in D^b (Y\times Z)$, let $b\circ a \in D^b (X\times Z)$ be 
$$b\circ a \cneq \dR p_{13\ast}(p_{12}^{\ast}(a)\dotimes p_{23}^{\ast}(b)).$$
It is easy to see $\Phi ^b _{Y\to Z}\circ \Phi _{X\to Y}^a \cong \Phi _{X\to Z}^{b\circ a}$. We have the following functor:
$$\pP \circ \colon D^b (X\times X)\ni a \longmapsto \pP \circ a \in D^b (X\times Y).$$
The above functor is an equivalence, since the functor 
$D^b (X\times Y) \ni b \mapsto \eE \circ b \in D^b (X\times X)$ gives a quasi-inverse. Here $\eE$ is a kernel of $\Phi ^{-1}$. 
Similarly we have an equivalence $\circ \pP \colon D^b (Y\times Y)\ni a \mapsto a\circ \pP \in D^b (X\times Y)$. Consider the 
following diagrams: ($\spadesuit$)
$$\begin{CD}
D^b (X\times X) @>{\pP \circ}>> D^b (X\times Y) \\
@A{\Delta _{X\ast}}AA   @AA{p_1 ^{\ast}(\ast)\dotimes \pP}A \\
D^b (X) @= D^b (X),
\end{CD} \qquad \quad
\begin{CD}
D^b (Y\times Y) @>{\circ\pP}>> D^b (X\times Y) \\
@A{\Delta _{Y\ast}}AA   @AA{p_2 ^{\ast}(\ast)\dotimes \pP}A \\
D^b (Y) @= D^b (Y).
\end{CD}$$
The above diagrams are 2-commutative. 
Let us check the left diagram commutes. Take $a\in D^b (X)$. Then 
\begin{align*}
\pP \circ (\Delta _{X\ast}a) &\cong \dR p_{13\ast}\left( p_{12}^{\ast}\Delta _{X\ast}a \dotimes p_{23}^{\ast}\pP \right) \\
&\cong \dR p_{13\ast}\left( (\Delta _X\times \id )_{\ast}p_1 ^{\ast}a \dotimes p_{23}^{\ast}\pP \right) \\
&\cong \dR p_{13\ast}(\Delta _X\times \id )_{\ast}\left( p_1 ^{\ast}a \dotimes (\Delta _X\times \id )^{\ast}p_{23}^{\ast}\pP \right) \\
&\cong p_1 ^{\ast}a \dotimes \pP.
\end{align*}
The second isomorphism follows from 
flat base change of the diagram below 
$$\begin{CD}
X\times Y @>{\Delta _X \times \id _Y}>> X\times X\times Y \\
@V{p_1}VV @VV{p_{12}}V \\
X @>{\Delta _X}>> X\times X, 
\end{CD}$$
and the third isomorphism is the projection formula. By the above commutative diagram, we have $\pP \circ \oO _{\Delta _X}\cong \pP$, 
$\oO _{\Delta _Y} \circ \pP \cong \pP$. Therefore we have the isomorphisms:
$$HH^{\ast}(X)\stackrel{\sim}{\longrightarrow} \Ext _{X\times Y}^{\ast}(\pP ,\pP) \stackrel{\sim}{\longrightarrow}HH^{\ast}(Y).$$
Since the equivalence $\Phi _{X\times X \to Y\times Y}^{p_{13}^{\ast}\pP \boxtimes p_{24}^{\ast}\eE}$ given in Theorem~\ref{hh} is nothing but the following
functor: 
$$\pP \circ (\ast) \circ \eE \colon D^b (X\times X)\longrightarrow D^b (Y\times Y),$$
the composition of the above isomorphisms is equal to $\phi$. $\quad$ q.e.d

\hspace{5mm}

Now let us take the exponential Atiyah class of $\pP$
$$\exp a(\pP) \colon \pP \longrightarrow \bigoplus _{i\ge 0}\pP \otimes \Omega _{X\times Y} ^i [i], $$
and take direct summands,
$$\exp a(\pP)_X \colon \pP \longrightarrow \bigoplus _{i\ge 0}\pP \otimes p_1 ^{\ast}\Omega _X ^i [i], \quad 
\exp a(\pP)_Y \colon \pP \longrightarrow \bigoplus _{i\ge 0}\pP \otimes p_2 ^{\ast}\Omega _Y ^i [i].$$
By the commutative diagram $(\spadesuit)$ in the proof of Lemma~\ref{co}, we have two morphisms
$$\exp (a)_X ^{+}\colon \oO _{\Delta _X}\longrightarrow \bigoplus _{i\ge 0}\Delta _{\ast}\Omega _X ^i [i], \quad 
\exp (a)_Y ^{+}\colon \oO _{\Delta _Y}\longrightarrow \bigoplus _{i\ge 0}\Delta _{\ast}\Omega _Y ^i [i],$$
such that $\pP \circ \exp (a)_X ^{+}=\exp a(\pP)_X$, $\exp (a)_Y ^{+}\circ \pP =\exp a(\pP)_Y$. 
We will investigate the relationship between $\exp (a)_X ^{+}$, $\exp (a)_Y ^{+}$, and the universal exponential Atiyah classes of $X$ and $Y$. Let 
$\sigma \colon X\times X \to X\times X$ be the involution $\sigma (x,x')=(x',x)$. 
\begin{lem}\label{exp} We have the following equalities:
$$\exp (a) _X ^{+} =\sigma _{\ast}\circ \exp (a)_X, \qquad \exp (a)_Y ^{+} =\exp (a)_Y.$$
\end{lem}
\textit{Proof}. 
We show $\exp (a)_X ^{+} =\sigma _{\ast}\circ \exp (a) _X$. 
 Let 
$$a _{X,i}^{+}\colon \oO _{\Delta _X}\longrightarrow \Delta _{\ast}\Omega _X ^i [i], \quad a(\pP )_{X,i}\colon \pP \longrightarrow
\pP \otimes p_1 ^{\ast}\Omega _X ^i [i]$$
be direct summands of $\exp (a)_X ^{+}$ and $\exp a(\pP) _X$ respectively. 
For $i=1$, we write $\ast _{1}=\ast$ for $\ast =a_X ^{+}$ or $a(\pP)_X$. 
We will show $a_{X,i}^{+}=\sigma _{\ast}a_{X,i}$. This is equivalent to $a(\pP)_{X,i}=\pP \circ (\sigma _{\ast}a_{X,i})$. 
First we treat the case of $i=1$.

Let $p_{ij}$ and $q_{ij}$ be projections from $X\times X\times Y$, $X\times Y \times X\times Y$ onto corresponding factors. Let 
\begin{align*} \Delta _X \times \id  &\colon X\times Y \hookrightarrow X\times X\times Y , \\
\id \times \Delta _Y &\colon X\times X\times Y \hookrightarrow X\times Y\times X\times Y 
\end{align*}
be $(\Delta _X \times \id)(x,y)=(x,x,y)$, $(\id \times \Delta _Y)(x,x',y)=(x,y,x',y)$. 
Let $I_{\Delta _{X ^{(2)},Y}}$ be the kernel of the composition
$$\oO _{X\times Y\times X\times Y}\longrightarrow (\id \times \Delta _Y)_{\ast}\oO _{X\times X\times Y}
\longrightarrow (\id \times \Delta _Y)_{\ast}\oO _{X\times X\times Y}/p_{12}^{\ast}I_{\Delta _X}^2.$$
Then we have a morphism of distinguished triangles, (in fact morphism of exact sequence)
$$\begin{CD}
\oO _{X\times Y\times X\times Y}/I_{\Delta_{X\times Y}}^2 @>>> \oO _{\Delta_{X\times Y}} @>>> \Delta _{X\times Y\ast}\Omega _{X\times Y}[1]\\
 @VVV @| @VVV \\
\oO _{X\times Y\times X\times Y}/I_{\Delta _{X ^{(2)},Y}} @>>> \oO _{\Delta_{X\times Y}} 
@>>> \Delta _{X\times Y\ast}p_1 ^{\ast}\Omega _X [1].
\end{CD} \quad (\diamondsuit)$$
Note that since 
\begin{align*} \Delta _{X\times Y\ast}p_1 ^{\ast}\Omega _X &\cong (\id \times \Delta _Y)_{\ast}(\Delta _X \times \id)_{\ast}p_1 ^{\ast}\Omega _X \\
&\cong (\id \times \Delta _Y)_{\ast}p_{12}^{\ast}\Delta _{X\ast}\Omega _X , \\
\end{align*}
and 
$$\oO _{X\times Y\times X\times Y}/I_{\Delta _{X ^{(2)},Y}} \cong
 (\id \times \Delta _Y)_{\ast}p_{12}^{\ast}\oO _{X\times X}/I_{\Delta _X}^2 ,
$$
the bottom sequence of $(\diamondsuit)$ is obtained by applying $(\id \times \Delta _Y)_{\ast} p_{12}^{\ast}$ to the 
distinguished triangle,
$$\Delta _{X\ast}\Omega _X \longrightarrow \oO _{X\times X}/I_{\Delta _X}^2 \longrightarrow \oO _{\Delta _X} 
\stackrel{a_X}{\longrightarrow} \Delta _{X\ast}\Omega _X [1].$$
Let $\widetilde{\Phi}\colon D^b(X\times Y\times X\times Y)\to D^b(X\times Y)$ be the functor $\widetilde{\Phi}\cneq \dR q_{34\ast}(\ast \dotimes q_{12}^{\ast}\pP)$. Then we have the isomorphisms of functors,

\begin{align*}& \widetilde{\Phi}\circ (\id \times \Delta _Y)_{\ast}\circ p_{12}^{\ast}(\ast) \\
&= \dR q_{34\ast}\left( (\id \times \Delta _Y)_{\ast}p_{12}^{\ast}(\ast)\dotimes q_{12}^{\ast}\pP \right) \\
&\cong \dR q_{34\ast}(\id \times \Delta _Y)_{\ast}\left( p_{12}^{\ast}(\ast)\dotimes (\id \times \Delta _Y)^{\ast}q_{12}^{\ast}\pP \right) \\
&\cong \dR p_{23\ast}(p_{12}^{\ast}(\ast)\dotimes p_{13}^{\ast}\pP) \\
&\cong \dR p_{23\ast}(\sigma \times \id)_{\ast}\left( (\sigma \times \id)^{\ast}p_{12}^{\ast}(\ast)\dotimes (\sigma \times \id)^{\ast}p_{13}^{\ast}\pP \right) \\
&\cong \dR p_{13\ast}(p_{12}^{\ast} \sigma _{\ast}(\ast)\dotimes p_{13}^{\ast}\pP) \\
&= \pP \circ \sigma _{\ast}(\ast).
\end{align*}
Therefore if we apply $\widetilde{\Phi}$ to the diagram $(\diamondsuit)$, we obtain the morphism of distinguished triangles,
$$\begin{CD} \widetilde{\pP} @>>> \pP @>{a(\pP)}>> \pP \otimes \Omega _{X\times Y}[1] \\
 @VVV @VVV @VVV \\
 \pP \circ \left( \sigma _{\ast}\oO _{X\times X}/I_{\Delta _X}^2 \right) 
@>>> \pP \circ (\oO _{\Delta _X}) @>{\pP \circ (\sigma _{\ast}a_X)}>> \pP \circ (\Delta _{X\ast}\Omega _X [1]).
\end{CD}$$
Here $\widetilde{\pP}\cneq \widetilde{\Phi}\left( \oO _{X\times Y\times X\times Y}/I_{\Delta _{X\times Y}}^2 \right)$. Since the morphism
$\pP \to \pP \circ (\oO _{\Delta _X})$, $\pP \otimes \Omega _{X\times Y} \to \pP \circ (\Delta _{X\ast}\Omega _X)$ of the 
above diagram are equal to $\id _{\pP}$, and direct summand $\pP \otimes \Omega _{X\times Y} \to \pP \otimes p_{1}^{\ast}\Omega _X$ under the isomorphism
$\pP \circ \Delta _{\ast} \cong \pP \dotimes p_1 ^{\ast}(\ast)$ of the 
diagram in Lemma~\ref{co}, we can conclude $a(\pP)_X =\pP \circ (\sigma _{\ast}a_X)$.

Secondly we show $a(\pP)_{X,i}=\pP \circ (\sigma _{\ast}a_{X,i})$ for all $i$. Since  
\begin{align*}
\pP \circ \sigma_{\ast}(a_{X}^{}\otimes p_2 ^{\ast}\Omega _X ^{\otimes i})&=
\pP \circ (\sigma _{\ast}a_X \otimes p_{1}^{\ast}\Omega _X ^{\otimes i}) \\
&= \dR p_{13\ast}(p_{12}^{\ast}\sigma_{\ast}a_{X}^{}\otimes p_{12}^{\ast}p_1 ^{\ast}\Omega _X ^{\otimes i}\dotimes p_{23}^{\ast}\pP )\\
&= \dR p_{13\ast}(p_{12}^{\ast}\sigma_{\ast}a_{X}^{}\otimes p_{13}^{\ast}p_1 ^{\ast}\Omega _X ^{\otimes i}\dotimes p_{23}^{\ast}\pP)\\
&= \dR p_{13\ast}(p_{12}^{\ast}\sigma _{\ast}a_{X}^{}\dotimes p_{23}^{\ast}\pP)\otimes p_1 ^{\ast}\Omega _X ^{\otimes i}\\
&=(\pP \circ \sigma _{\ast}a_{X}^{})\otimes p_1 ^{\ast}\Omega _X ^{\otimes i} \\
&= a(\pP)_{X}\otimes p_1 ^{\ast}\Omega _X ^{\otimes i},
\end{align*}
we have $a(\pP)_{X,i}=\pP \circ (\sigma _{\ast}a_{X,i})$. $\quad$ q.e.d

\hspace{5mm}

Using the above proposition, we can find the relationship between HKR-isomorphism, the isomorphism $HH^{\ast}(X)\to \Ext _{X\times Y}^{\ast}(\pP, \pP)$
of Lemma~\ref{co} and the exponential Atiyah-classes. In fact we have the following lemma:
\begin{lem}\label{coo} The following diagrams commute: 
\begin{align*}\begin{CD}HT^{\ast}(X\times Y) @>{\times \exp a(\pP)}>> \Ext ^{\ast}_{X\times Y}(\pP ,\pP) \\
@A{p_1 ^{\ast}}AA  @AA{\pP \circ }A \\
HT^{\ast}(X) @>>{\sigma _{\ast}I_X}> HH^{\ast}(X),
\end{CD} \\
\begin{CD}HT^{\ast}(X\times Y) @>{\times \exp a(\pP)}>> \Ext ^{\ast}_{X\times Y}(\pP ,\pP) \\
@A{p_2 ^{\ast}}AA  @AA{\circ \pP}A \\
HT^{\ast}(Y) @>>{I_Y}> HH^{\ast}(Y).
\end{CD}\end{align*}
Here $\times \exp a(\pP)$ means multiplying by $\exp a(\pP)$ and taking $\Ext ^{\ast}(\pP, \pP)$-component. 
\end{lem}
\textit{Proof}. We show that the top diagram commutes. Take $u\in H^p (X, \wedge ^q T_X)$. By Lemma~\ref{exp}, we have $\sigma _{\ast} a_{X,q}=a_{X,q}^{+}$. So by Proposition~\ref{ati}, 
$\sigma _{\ast}I_X (u)$ is the composition:
$$\oO _{\Delta}\stackrel{a_{X,q}^{+}}{\longrightarrow}\Delta _{\ast}\Omega _X ^q [q] \stackrel{\Delta _{\ast}u}{\longrightarrow}
\Delta _{\ast}\oO _X [p+q].$$
Therefore $\pP \circ \sigma _{\ast}I_X (u)$ is the composition
$$\pP \stackrel{a(\pP)_{X,q}}{\longrightarrow}\pP \otimes p_1 ^{\ast}\Omega _X ^q [q] 
\stackrel{\times p_1 ^{\ast}u}{\longrightarrow}\pP [p+q].$$
But this is equal to the composition
$$\pP \stackrel{a(\pP )_q}{\longrightarrow}\pP \otimes \Omega _{X\times Y}^q [q] \stackrel{\times p_1 ^{\ast}u}{\longrightarrow}
\pP [p+q].$$
Therefore the diagram commutes. $\quad$ q.e.d

\section{Proof of the main theorem}
In this section we will prove Theorem~\ref{mt}.
 Let $X,Y$ and $\Phi$, $\pP$ be as in the previous sections.
We want to extend $\pP$ to $\pP ^{\dag}\in D^b _{\perf}(X\times Y, -p_1 ^{\ast}\check{u}+p_2 ^{\ast}v)$. 
For this purpose we have to investigate the relationship between $u$, $v$, and the exponential Atiyah-class of $\pP$.  
For $u\in H^p (X,\wedge ^q T_X)$, let $\check{u} \cneq (-1)^q u$, and extend the operation to $HT^{\ast}(X)$ linearly.  
Then it is clear that $\sigma _{\ast}I_X (u)=I_X (\check{u})$. 
Take $u\in HT^{\ast}(X)$ and $v=\phi _{T}(u)$. By Lemma~\ref{coo} and the above remark, we have 
\begin{align*}(-p_1 ^{\ast}\check{u}+p_2 ^{\ast}v) \cdot \exp a(\pP) &= -\pP \circ \sigma _{\ast}I_X (\check{u}) +I_Y (v)\circ \pP \\
&= -\pP \circ I_X (u)+ \left( \phi \circ I_X (u) \right) \circ \pP \\
&= -\pP \circ I_X (u)+ \pP \circ I_X (u) \circ \eE \circ \pP \\
&=0, 
\end{align*}
in $\Ext _{X\times Y}^2 (\pP, \pP)$. Therefore to extend $\pP$ to $\pP ^{\dag}$, it suffices to show the following proposition.
\begin{prop}\label{ext}
Take $\pP \in D^b (X)$ and $u\in HT^2 (X)$. Assume that $u\cdot \exp a(\pP)=0$ 
in $\Ext _X ^2 (\pP, \pP)$. Then there exists an object $\pP ^{\dag}\in D^b _{\perf}(X,u)$ such that 
$\dL i^{\ast}\pP ^{\dag} \cong \pP$. 
\end{prop}
\textit{Proof}.
Let $\pP ^{\bullet}$ be a complex of locally free sheaves on $X$, which represents $\pP$. 
Since $\pP ^n$ is locally free, we have 
$$\widetilde{\pP}^i \cneq \dR p_{2\ast} \left( p_1 ^{\ast}\pP ^n \dotimes \oO _{X\times X}/I_{\Delta _X}^2 \right) 
=p_{2\ast}\left( p_1 ^{\ast}\pP ^n \otimes \oO _{X\times X}/I_{\Delta}^2 \right),$$
and the distinguished triangle 
$$\pP ^{}\otimes \Omega _X \longrightarrow \dR p_{2\ast}  \left( p_1 ^{\ast}\pP ^{} \dotimes \oO _{X\times X}/I_{\Delta _X}^2 \right)  
\longrightarrow \pP ^{} \longrightarrow \pP ^{}\otimes \Omega _X [1]$$
is represented by the exact sequence of complexes, 
$$0\longrightarrow \pP ^{\bullet}\otimes \Omega _X \longrightarrow \widetilde{\pP}^{\bullet}\stackrel{\psi ^{\bullet}}{\longrightarrow}
 \pP ^{\bullet} \longrightarrow 0.$$
But $\psi ^n \colon 
\widetilde{\pP} ^n \to \pP ^n$ has a $\mathbb{C}$-linear section $\lambda ^n 
 \colon \pP ^n \to \widetilde{\pP }^n$,
$$\lambda ^n (U_i)\colon \pP ^n (U_i)\ni x \longmapsto x\otimes 1 \in \widetilde{\pP}^n (U_i) =\pP ^n (U_i)\otimes _{\oO _{U_i}} \oO _{U_i \times U_i}/I_{\Delta}^2 $$
and $\lambda ^{\bullet}\colon \pP ^{\bullet}\to \widetilde{\pP}^{\bullet}$ gives a $\mathbb{C}$-linear splitting of $\psi ^{\bullet}\colon 
\widetilde{\pP}^{\bullet}\to \pP ^{\bullet}$.
Therefore the 
Atiyah class $a(\pP)$ becomes the zero map after  
applying the forgetful functor $D^b (X)\to D^b (\Mod (X,\mathbb{C}))$. 
Here $\Mod (X, \mathbb{C})$ is a category of sheaves of 
$\mathbb{C}$-vector spaces 
on $X$.

On the other hand, in the derived category of quasi-coherent sheaves, 
 the Atiyah class is represented by some morphism of complexes
 of quasi-coherent sheaves, denoted by the same symbol $a(\pP)$:
$$a(\pP) \colon \pP ^{\bullet}\longrightarrow T \cC ^{\bullet }(\mathfrak{U}, \pP ^{\bullet}\otimes \Omega _X). $$
Here $T$ is a translation functor 
$T(X^{\bullet})=X^{\bullet +1}$, $T(d_X)=-d_X$. 
By the above remark, $a(\pP)$ is homotopic to zero as 
complexes of $\mathbb{C}$-vector spaces, 
and we are now going to construct a homotopy. 
Let us choose connections 
$$\nabla _i ^{(n)} \colon \pP ^n |_{U_i}  \longrightarrow \pP ^n |_{U_i} \otimes \Omega _X $$
on $U_i$ for all $i$. 
Then it is easy to check that a homotopy between $a(\pP)$ and zero as 
morphisms of complexes of $\mathbb{C}$-vector spaces
 is given by 
$\nabla \colon \pP ^{\bullet}\to \cC ^{\bullet}(\mathfrak{U}, \pP ^{\bullet}\otimes 
\Omega _X)$, defined as follows:
$$\nabla ^n \colon \pP ^n \ni x \longmapsto \{ \nabla _i ^{(n)}(x) \}_i \in \cC ^0 (\mathfrak{U}, \pP ^n \otimes \Omega _X)\subset 
\cC ^n (\mathfrak{U}, \pP ^{\bullet}\otimes \Omega _X).$$
Namely $a(\pP)= \nabla \circ d_{\pP} +T(d_{\cC})\circ \nabla$. 
(cf.~\cite{Hu}). 
Also $a(\pP )_2$ is represented by a morphism of complexes of quasi-coherent 
sheaves, 
$$a(\pP)_2 \colon \pP ^{\bullet}\longrightarrow T^2 \cC ^{\bullet}(\mathfrak{U}, \pP ^{\bullet}\otimes \Omega _X ^2)=
\cC ^{\bullet +2}(\mathfrak{U}, \pP ^{\bullet}\otimes \Omega _X ^2)$$
which is homotopic to zero as complexes of $\mathbb{C}$-vector spaces.
 In fact
we can calculate $a(\pP)_2$ as follows: 
\begin{align*}a(\pP)_2 &= \epsilon \circ T(a(\pP)\otimes 1)\circ a(\pP) \\
&= -\epsilon(\nabla \circ d_{\cC}\circ \nabla)\circ d_{\pP}-d_{\cC}\circ \epsilon (\nabla \circ d_{\cC} \circ \nabla).
\end{align*}
Hence the homotopy is given by 
$$-\epsilon(\nabla \circ d_{\cC}\circ \nabla)\colon \pP ^{\bullet}
\lr T\cC ^{\bullet}(\mathfrak{U}, \pP ^{\bullet}\otimes \Omega _X ^2).$$ 
Here $d_{\pP}$ and $d_{\cC}$ are differentials of $\pP ^{\bullet}$ and $\cC ^{\bullet}(\mathfrak{U}, \pP ^{\bullet})$ respectively. 
Therefore 
if we take a $\check{\mbox{C}}$ech representative of $\beta$ 
and consider the morphism of complexes of 
quasi-coherent sheaves
 $$\beta \cdot a(\pP) +\gamma \cdot a(\pP)_2 \colon  \pP ^{\bullet}\longrightarrow \cC ^{\bullet +2}(\mathfrak{U}, \pP ^{\bullet}),$$
then this is homotopic to zero as 
morphisms of complexes of $\mathbb{C}$-vector spaces. The homotopy is given by 
$$\nabla ^{\dag}\cneq \beta \circ \nabla -\gamma \circ \epsilon (\nabla \circ d_{\cC}\circ \nabla)\colon \pP ^{\bullet} \longrightarrow 
T\cC ^{\bullet}(\mathfrak{U}, \pP ^{\bullet}).$$
By the assumption, 
$$\alpha \otimes \id _{\pP} + \beta \cdot a(\pP)  +\gamma \cdot a(\pP)_2  \colon \pP ^{\bullet}\longrightarrow 
\cC ^{\bullet +2}(\mathfrak{U}, \pP ^{\bullet})$$
is homotopic to zero as a map of complexes of quasi-coherent sheaves, and let 
$h ^{\bullet} \colon \pP ^{\bullet}\longrightarrow T\cC ^{\bullet }(\mathfrak{U}, \pP ^{\bullet})$
be such a homotopy. Note that $h^n$ is a $\oO _X$-module homomorphism. Combining these, we can conclude $\alpha \otimes \id _{\pP}$ is homotopic to 
zero as complexes of $\mathbb{C}$-vector spaces and the 
homotopy is given by $h^{\dag}\cneq h-\nabla ^{\dag}$.

 Now we are going to 
  construct the complex $(\pP ^{\dag})^{\bullet}$ whose terms are
 objects in $\QCoh (X,u)$ by using $h^{\dag}$. First define $(\pP ^{\dag})^n _i$ to be
 $$(\pP ^{\dag})^n _i \cneq \pP ^n |_{U_i} \oplus \cC ^n (\mathfrak{U}, \pP ^{\bullet})|_{U_i}.$$
 We introduce a left $\oO ^{(\beta, \gamma)}_{U_i}$-module structure on $(\pP ^{\dag})^n _i$. For $a\in \oO _{U_i}$, let $\gamma _a \in T_{U_i}$ be 
 a differential operator $\gamma _a \cneq \gamma (a, \ast)$. Then for 
$$(a,b) \in \oO _{U_i}\oplus \cC ^0 (\mathfrak{U}, \oO _X)|_{U_i}, \quad 
(x,y)  \in \pP ^n |_{U_i} \oplus \cC ^n (\mathfrak{U}, \pP ^{\bullet})|_{U_i},
$$
define $(a,b)\ast _{\gamma}(x,y)$ to be 
$$(a,b)\ast _{\gamma}(x,y)\cneq (ax, bx +ay + \{ \gamma _a \circ \nabla _{i_0}^{(n)}(x|_{U_{i_0}})\}_{i_0} )\in (\pP ^{\dag})^n _i .$$ 
Here $ \{ \gamma _a \circ \nabla _{i_0}^{(n)}(x|_{U_{i_0}})\}_{i_0} \in \cC ^0 (\mathfrak{U}, \pP ^n)|_{U_i}$. We have to 
check the following:
\begin{lem}\label{act}
$\ast _{\gamma}$ defines the left action of $\oO _{U_i} ^{(\beta, \gamma)}\subset \oO _X |_{U_i}\oplus \cC ^0 (\mathfrak{U}, \oO _X)|_{U_i}$ on 
$(\pP ^{\dag})^n _i$.  
\end{lem}
\textit{Proof}. Take $(a,b), (c,d) \in  \oO _X ^{(\beta, \gamma)}(U_i)$ and $(x,y)\in (\pP ^{\dag})^n _i$. 
Then 
\begin{align*}
& (a,b)\ast _{\gamma} \{ (c,d)\ast _{\gamma} (x,y) \} \\
& =(a,b)\ast _{\gamma} (cx, cy+dx+ \{ \gamma _c \circ \nabla _{i_0}^{(n)} (x|_{U_{i_0}}) \} _{i_0}) \\
& = (acx, acy+adx+bcx+\{ a\gamma _c \circ \nabla _{i_0}^{(n)}(x|_{U_{i_0}}) +\gamma _a \circ \nabla _{i_0} ^{(n)} (cx |_{U_{i_0}}) \} _{i_0}),
\end{align*} 
and 
\begin{align*}
&\{ (a,b)\ast _{\gamma}(c,d)\} \ast _{\gamma} (x,y) \\
& =(ac, \gamma (a,c) +ad+bc)\ast _{\gamma} (x,y)\\
&= (acx, acy + adx +bcx + \gamma (a,c) \cdot x+\{ \gamma _{ac}\circ \nabla _{i_0}^{(n)} (x|_{U_{i_0}})\}_{i_0} )
\end{align*}
Since $\nabla _{i_0}^{(n)}$ is a connection, we have 
\begin{align*}
& a\gamma _c \circ \nabla _{i_0}^{(n)} (x|_{U_{i_0}})+\gamma _{a}\circ \nabla _{i_0}^{(n)}(cx|_{U_{i_0}}) \\
&= a\gamma _c \circ \nabla _{i_0}^{(n)} (x|_{U_{i_0}})+\gamma _a \circ \{ dc \otimes (x|_{U_{i_0}})+c\cdot \nabla _{i_0}^{(n)}(x|_{U_{i_0}}) \} \\
&= \gamma (a,c)\cdot (x|_{U_{i_0}})+(c\gamma _a +a \gamma _c)\circ \nabla _{i_0}^{(n)}(x|_{U_{i_0}}) \\
&= \gamma (a,c)\cdot (x|_{U_{i_0}})+\gamma _{ac}\circ \nabla _{i_0}^{(n)}(x|_{U_{i_0}})
\end{align*} 
Therefore the lemma follows. $\quad$ q.e.d

\hspace{5mm}

By Lemma~\ref{act}, we have obtained the object, 
$$(\pP ^{\dag})^n _i \in \Mod (\oO _{U_i}^{(\gamma)}).$$
If we regard $\pP ^n |_{U_i}$ and $\cC ^n (\mathfrak{U}, \pP ^{\bullet})
|_{U_i}$ as $\oO_{U_i}^{(\gamma)}$-modules by the surjection, 
$\oO _{U_i}^{(\gamma)}\twoheadrightarrow \oO _{U_i}$, then we 
have the exact sequence in $\Mod (\oO _{U_i}^{(\gamma)})$, 
$$0\longrightarrow \cC ^n (\mathfrak{U}, \pP ^{\bullet})
|_{U_i}\longrightarrow (\pP ^{\dag})^n _i \longrightarrow 
\pP ^n |_{U_i}\longrightarrow 0.$$
Since $\pP ^n |_{U_i}$, $\cC ^n (\mathfrak{U}, \pP ^{\bullet})
|_{U_i}$ are objects in $\QCoh (\oO _{U_i}^{(\gamma)})$, we have 
$$(\pP ^{\dag})^n _i \in \QCoh (\oO _{U_i}^{(\gamma)}),$$
by Lemma~\ref{coh}(3) and Corollary~\ref{ne}.
Next define $\phi _{i_0 i_1}^n \colon (\pP ^{\dag})^n _{i_0}|_{U_{i_0 i_1}} \to  (\pP ^{\dag})^n _{i_1}|_{U_{i_0 i_1}}$ to be
$$\phi _{i_0 i_1}^n (x,y)\cneq (x, -\{ \alpha _{i_0 i_1 j}\cdot x \}_j +y).$$
Here $-\{ \alpha _{i_0 i_1 j}\cdot x \}_j \in \cC ^0 (\mathfrak{U}, \pP ^n)|_{U_{i_0 i_1}}$. Then $\phi _{i_0 i_1} ^n$ is clearly 
$\oO _{U_{i_0 i_1}} ^{(\beta, \gamma)}$-module homomorphism, and the cocycle condition of $\alpha$ implies the following:
$$(\pP ^{\dag})^n \cneq ((\pP ^{\dag})^n _i, \phi _{i_0 i_1}^n ) \in \QCoh (X, u).$$
Now we will construct a differential $d^n \colon (\pP ^{\dag})^n \to (\pP ^{\dag})^{n+1}$. On $U_i$, we define $d_i ^n$ as 
$$d_i ^n (x,y)\cneq (d_{\pP}x, d_{\cC}y +h^{\dag}(x)-\{\alpha _{ikl} \cdot x \}_{kl} )
\in \pP ^{n+1}|_{U_i} \oplus \cC ^{n+1} (\mathfrak{U}, \pP ^{\bullet} )|_{U_i}, $$ 
for $(x,y) \in \pP ^{n}|_{U_i} \oplus \cC ^{n} (\mathfrak{U}, \pP ^{\bullet} )|_{U_i}$.  
Here $ \{\alpha _{ikl} \cdot x \}_{kl} \in \cC ^1 (\mathfrak{U}, \pP ^n)|_{U_i}$. Then 
\begin{align*}
& d_i ^{n+1}\circ d_i ^n (x,y)\\
&= (0, d_{\cC}(d_{\cC}y+h^{\dag}(x)-\{ \alpha _{ikl}\cdot x \} _{kl})+h^{\dag}(d_{\pP}x)-\{\alpha _{ikl}\cdot d_{\pP}x \}_{kl}) \\
&= (0, d_{\cC}h^{\dag}(x)+h^{\dag}d_{\pP}(x)-\{\alpha _{i_0 i_1 i_2}\cdot x\}_{i_0 i_1 i_2})\\
&=(0,0).
\end{align*}
The second equality comes from the cocycle condition of $\alpha$. We can check $\phi ^{n+1}_{i_0 i_1}\circ d_{i_0}^n =
d_{i_1}^n \circ \phi _{i_0 i_1}^n$ similarly.  We have to check the following:
\begin{lem}$d_i ^n$ is $\oO _{U_i} ^{(\beta, \gamma)}$-module homomorphism. \end{lem}
\textit{Proof}. Take $(a,b)\in \oO _{U_i} ^{(\beta, \gamma)}\subset \oO _{U_i}\oplus \cC ^0 (\mathfrak{U},\oO _X)|_{U_i}$, i.e.
$\delta b=\{ \beta _{i_0 i_1}(a)\}_{i_0 i_1}$, and $(x,y)\in \pP ^n|_{U_i} \oplus \cC ^n (\mathfrak{U}, \pP ^{\bullet})|_{U_i}$. 
Then 
\begin{align*}
& (a,b)\ast _{\gamma}d_i ^n(x,y) \\
&= (a,b)\ast _{\gamma}(d_{\pP}x, d_{\cC}y +h^{\dag}(x)-\{ \alpha _{ijk}\cdot x \} _{jk})\\
&= (ad_{\pP}x, ad_{\cC}y +ah^{\dag}(x)-a\{\alpha _{ijk}\cdot x \}_{jk}+bd_{\pP}x+\{ \gamma _a \circ \nabla _{i_0}^{(n+1)}(d_{\pP}x)\} _{i_0}),
\end{align*}
and 
\begin{align*}
& d_i ^n \{ (a,b)\ast _{\gamma}(x,y)\} \\ 
&= d_i ^n(ax, ay +bx +\{ \gamma _a \circ \nabla _{i_0}^{(n)}(x)\} _{i_0}) \\
&= (d_{\pP}(ax), d_{\cC}(ay+bx+\{ \gamma _a \circ \nabla _{i_0}^{(n)}(x)\} _{i_0})+ h^{\dag}(ax)-\{ \alpha _{ijk}\cdot ax\} _{jk}).
\end{align*} 
Therefore it suffices to check the following:
$$-a\nabla ^{\dag}(x)+\{ \gamma _a \circ \nabla _{i_0}^{(n+1)}(d_{\pP}x) \} _{i_0}
=\delta (bx)+d_{\cC}\{ \gamma _a \circ \nabla _{i_0}^{(n)}(x)\} _{i_0}-\nabla ^{\dag}(ax).$$
We calculate $\nabla ^{\dag}(ax)-a\nabla ^{\dag}(x)$. Since  
$$\nabla (ax)-a \nabla (x)=da \otimes x,$$
and 
\begin{align*}
& \nabla \circ d_{\cC} \circ \nabla (ax)-a\nabla \circ d_{\cC}\circ \nabla (x) \\
&= \nabla \circ d_{\cC}\circ (da \otimes x +a\nabla (x))-a\nabla \circ d_{\cC}\circ \nabla (x) \\
&= \nabla \circ (da\otimes d_{\pP}x+ ad_{\cC}\circ \nabla (x))-a\nabla \circ d_{\cC}\circ \nabla (x)\\
&= da \otimes \nabla \circ d_{\pP}x +d_{\cC}\circ \nabla (x)\otimes da, 
\end{align*}
we have 
\begin{align*}
& \nabla ^{\dag}(ax)-a\nabla ^{\dag}(x) \\
&= \{ \beta _{i_0 i_1}(a)\} _{i_0 i_1}\cdot x -\gamma _a \circ \nabla \circ d_{\pP}(x)+\gamma _a \circ d_{\cC}\circ \nabla (x)\\
&= \delta (b)x-\gamma _a \circ \nabla \circ d_{\pP}(x)+\gamma _a \circ d_{\cC}\circ \nabla (x).
\end{align*}
So the lemma follows. $\quad$ q.e.d

\hspace{5mm}

We have constructed an unbounded complex of $\QCoh (X,u)$:
$$\pP ^{\dag}\cneq 
\cdots \lr (\pP ^{\dag})^n \stackrel{d^n}{\lr} (\pP ^{\dag})^{n+1} \lr 
\cdots,$$
with $d^n |_{U_i}=d_i ^n$. 
 The next lemma finishes the proof of Proposition~\ref{ext}. $\quad$ q.e.d 

\begin{lem}
$\pP ^{\dag}$ is locally quasi-isomorphic to a bounded complex of free $\oO _X ^{(\beta, \gamma)}$ -modules of finite rank, and 
$\dL i^{\ast}\pP ^{\dag}\cong \pP$. \end{lem}
\textit{Proof}. Let 
$$p_i ^n \colon \cC ^n (\mathfrak{U}, \pP ^{\bullet})|_{U_i} \longrightarrow \cC ^0 (\mathfrak{U}, \pP ^n )|_{U_i} \longrightarrow \pP ^n |_{U_i}$$
be a projection and $\widetilde{h}^n _i$ be the composition,
$$\widetilde{h}^n _i \cneq p_i ^{n+1}\circ (h^{\dag})^n _i \colon \pP ^n |_{U_i} \longrightarrow \cC ^{n+1}(\mathfrak{U}, \pP ^{\bullet})|_{U_i} 
\longrightarrow \pP ^{n+1}|_{U_i}.$$
Let $\widetilde{\pP}^n _i \cneq \pP ^n |_{U_i}[\varepsilon]/(\varepsilon ^2)$ be a free left $\oO _{U_i}^{(\gamma)}$-module, the left action 
given by for $a+b\varepsilon \in \oO _{U_i}^{(\gamma)}$, $x+y\varepsilon \in \widetilde{\pP}^n _i $, 
$$(a+b\varepsilon )\ast _{\gamma} (x+y\varepsilon)\cneq ax + (ay+bx+\gamma _a \circ \nabla _i ^{(n)}(x))\varepsilon.$$
Then define the complex $\widetilde{\pP}^{\bullet}_i $ whose differential is given by  
$$\widetilde{\pP}^n _i \ni x+y\varepsilon \longmapsto d_{\pP}(x)+(d_{\cC}y+\widetilde{h}^n _i (x))\varepsilon \in \widetilde{\pP}^{n+1} _i.$$
We will show the natural map
$$(\pP ^{\dag})^n _i \ni (x,y)\longmapsto x+p^n _i (y)\varepsilon \in \widetilde{\pP}^n _i $$
gives a quasi-isomorphism
between $\pP ^{\dag}|_{U_i}$ and $\widetilde{\pP}^{\bullet}_i$. 
 It is clear that the above map is a morphism of complexes of $\oO ^{(\beta ,\gamma)}_{U_i}$-modules. Note that 
$p_i ^{\bullet}$ gives a splitting of the $\check{\mbox{C}}$ech resolution on $U_i$, so we have the decomposition,
$$\cC ^{\bullet}(\mathfrak{U}, \pP ^{\bullet})|_{U_i}\cong \pP ^{\bullet}|_{U_i}\oplus \qQ _i ^{\bullet},$$
for some complex $\qQ ^{\bullet}_i$ with $H^{\bullet}(\qQ _i ^{\bullet})=0$. We have the following diagram:
$$\begin{CD} @. T\qQ ^{\bullet}_i @. \\
@. @VVV @. \\
\pP ^{\bullet}|_{U_i} @>{h^{\dag}-\{ \alpha _{ijk} \} _{jk}}>> T\cC ^{\bullet}(\mathfrak{U}, \pP ^{\bullet})|_{U_i} @>>> T(\pP ^{\dag})^{\bullet}_i \\
@| @VV{Tp_i ^{\bullet}}V @. \\
\pP ^{\bullet}|_{U_i} @>{\widetilde{h}}>> T\pP ^{\bullet}|_{U_i} @>>> T\widetilde{\pP}^{\bullet}_i .
\end{CD}$$ 
Therefore we obtain the distinguished triangle
$$\qQ ^{\bullet}_i \longrightarrow (\pP ^{\dag})^{\bullet}_i \longrightarrow \widetilde{\pP}^{\bullet}_i \longrightarrow T\qQ ^{\bullet}_i .$$
Since $\qQ ^{\bullet}_i$ is acyclic, the first part of the lemma follows. 
For the second part, we have a morphism of complexes $(\pP ^{\dag})^{\bullet}\to i_{\ast}\pP ^{\bullet}$ by construction. By taking adjoint, we have 
a morphism $\dL i^{\ast}\pP ^{\dag} \to \pP$ in $D^b (X)$. This morphism is quasi-isomorphic on $U_i$, hence quasi-isomorphic. 
$\quad$ q.e.d

\hspace{5mm}

Now let us return to the situation of the first part of this section. 
By Proposition~\ref{ext}, we obtain the object $\pP ^{\dag} \in D^b _{\perf} (X\times Y, -p_1 ^{\ast}\check{u}+p_2 ^{\ast}v)$. 
Therefore we can construct a functor $\Phi ^{\dag}\colon D^b (X,u)\to D^b (Y,v)$. Next we will show $\Phi ^{\dag}$ fits some commutative diagram. 
\begin{lem} \label{2com}
The following diagram is 2-commutative,
$$\begin{CD}
D^b (X) @>{i_{\ast}}>> D^b (X,{u}) @>{\dL i^{\ast}}>> D^{-}(X)\\
@V{\Phi}VV @VV{\Phi ^{\dag}}V @VV{\Phi ^{-}}V\\
D^b (Y) @>{i_{\ast}}>> D^b (Y, v)@>{\dL i^{\ast}}>> D^{-}(Y).
\end{CD}$$
\end{lem}
\textit{Proof}. To distinguish the notation, let 
\begin{align*}\dR p_{2\ast} ^{\dag} & \colon D^b (\Mod (p_2 ^{-1}\oO _Y ^{(\beta ',\gamma')},p_2 ^{\ast}\widetilde{\alpha}'))
\longrightarrow D^b (\Mod (\oO _Y ^{(\beta ',\gamma ')}),\widetilde{\alpha} '), \\
p_1 ^{\dag \ast}& \colon D^b (X,u)\longrightarrow D^b (X\times Y, p_1 ^{\ast}u +p_2 ^{\ast}(0,\beta ',\gamma ')),
\end{align*}
be derived push-forward and pull-back. 
Let us take $a\in D^b (X)$. Then 
\begin{align*}\Phi ^{\dag}\circ i_{\ast}(a) &= \dR p_{2\ast}^{\dag}(p_1 ^{\dag \ast}i_{\ast}a\dotimes \pP ^{\dag}) \\
&\cong \dR p_{2\ast}^{\dag}(i_{\ast}p_1 ^{\ast}a\dotimes \pP ^{\dag}) \\
&\cong \dR p_{2\ast}^{\dag}i_{\ast}(p_1 ^{\ast}a\dotimes \dL i^{\ast}\pP ^{\dag}) \\
&\cong i_{\ast}\dR p_{2\ast}(p_1 ^{\ast}a \dotimes \pP) \\
&\cong i_{\ast}\circ \Phi (a).
\end{align*}
The second isomorphism follows from flat base change, and the third from projection formula. These properties are verified in our case 
as in the commutative case. 
We have proved the left diagram commutes. The right diagram commutes similarly. $\quad$q.e.d

\hspace{5mm}

\textit{Proof of Theorem~\ref{mt}}.
It remains to show $\Phi ^{\dag}$ gives an equivalence. Take $a\in D^b (X,u)$ and $b\in D^{-}(X)$. Then we have 
\begin{align*} 
\Hom (\Phi^{\dag}(a),\Phi ^{\dag} i_{\ast}(b)) &\cong \Hom (\Phi ^{\dag}(a),i_{\ast}\Phi (b)) \\
&\cong \Hom (\dL i^{\ast}\Phi ^{\dag}(a),\Phi (b)) \\
&\cong \Hom (\Phi ^{-}\dL i^{\ast}a, \Phi (b)) \\
&\cong \Hom (\dL i^{\ast}a,b) \\
&\cong \Hom (a,i_{\ast}b).
\end{align*}
Therefore the map $\Hom (a,i_{\ast}b)\stackrel{\Phi ^{\dag}}{\longrightarrow}\Hom (\Phi ^{\dag}(a),\Phi ^{\dag}(i_{\ast}b))$ is an 
isomorphism. Next take $a,b \in D^b (X,u)$. Since we have the distinguished triangle,
$$i_{\ast}\dL i^{\ast}b \longrightarrow b \longrightarrow i_{\ast}\dL i^{\ast}b \longrightarrow i_{\ast}\dL i^{\ast}b[1],$$
we have the following morphism of exact sequences, ($b'\cneq \dL i^{\ast}b$)
$$\begin{CD}
 \Hom (a,i_{\ast}b') @>>> \Hom (a,b) @>>> \Hom (a,i_{\ast}b')  \\
 @VVV @VVV @VVV  \\
 \Hom (\Phi ^{\dag}(a),\Phi ^{\dag}(i_{\ast}b')) @>>> \Hom (\Phi ^{\dag}(a),\Phi ^{\dag}(b)) @>>> \Hom (\Phi ^{\dag}(a),\Phi ^{\dag}(i_{\ast}b'))   .
\end{CD}$$
Therefore the morphism $\Hom (a,b)\longrightarrow \Hom (\Phi ^{\dag}(a),\Phi ^{\dag}(b))$ is an isomorphism by 5-lemma. 
Now we have proved $\Phi ^{\dag}$ is fully-faithful. Finally we show $\Phi ^{\dag}$ is essentially surjective. 
Take ${\fF}\in D^b (Y,v)$. Again we have the distinguished triangle, 
$$i_{\ast}\dL i^{\ast}{\fF} \longrightarrow {\fF} \longrightarrow i_{\ast}\dL i^{\ast}{\fF} \stackrel{t_{\fF}}{\longrightarrow}i_{\ast}\dL i_{\ast}{\fF}[1].$$
Let ${\fF}'\cneq \dL i^{\ast}{\fF}$. Since we have 
\begin{align*}
i_{\ast}{\fF}' &\cong i_{\ast}\Phi \circ \Psi ({\fF}') \\
&\cong \Phi ^{\dag}\circ i_{\ast}\Psi ({\fF}'),
\end{align*}
the morphism $t_{\fF} \colon i_{\ast}{\fF}' \to i_{\ast}{\fF}'[1]$ is obtained 
by applying $\Phi ^{\dag}$ to some morphism,
$s_{\fF} \colon i_{\ast}\Psi ({\fF}')\to i_{\ast}\Psi ({\fF}')[1]$. Let $\gG\cneq \Cone (s_{\fF})$. Then ${\fF}$ is isomorphic to $\Phi ^{\dag}(\gG)$. 
It remains to show $\gG$ is bounded. Note that by the definition of 
$\Phi ^{\dag}$, there exists $N>0$ such that if $H^i (A)=0$ for $i\ge l$ and 
some $l$, then $H^i (\Phi ^{\dag}(A))=0$ for $i\ge l+N$. Let us 
take an intelligent truncation of $\gG$:
$$\tau _{\le l-1}\gG \longrightarrow \gG \longrightarrow \tau _{\ge l}\gG.$$
Then by the above remark, $H^i (\Phi ^{\dag}(\tau _{\le l-1}\gG))=0$ for 
$i\ge l+N$. Therefore $\Phi ^{\dag}(\tau _{\le l-1}\gG)\to \Phi ^{\dag}(\gG)
=\fF$ is zero-map for sufficiently small $l$. Since $\Phi ^{\dag}$ is fully-faithful, this implies $\tau _{\le l-1}\gG \to \gG$ is zero-map. Therefore 
$\tau _{\le l-1}\gG=0$. 
$\quad$ q.e.d

\section{Examples}
\subsection*{Abelian varieties}
We give an example in which $\phi _T$ does not preserve direct summands of 
$HT ^2 (X)$. 
Let $A$ be an Abelian variety, and $\hat{A}$ be its dual Abelian 
variety. Let $\uU \in \Pic (A\times \hat{A})$ be 
the Poincare line bundle. Then the functor
$$\Phi _{\hat{A}\to A}^{\uU}\colon D(\hat{A}) \lt D(A)$$
gives an equivalence. (cf.~\cite{Mu1}). In this particular 
example, $\phi _T$ takes some $\alpha \in H^2 (\oO _{\hat{A}})$ 
to $\gamma \in H^0 (\wedge ^2 T_A)$. Hence $\Phi ^{\dag}$ 
give equivalences between gerby deformations and non-commutative 
deformations of Abelian varieties first orderly. 
This phenomenon
 has been extended to infinite order deformations in~\cite{OJT}. 

\subsection*{Birational geometry}
In this example, we discuss the situation in which 
 $\phi _T$ preserves some direct summands of $HT^2 (X)$. 
 This example comes from the equivalences under some birational 
 transforms, e.g. flops. 
Recently the relationship between derived categories and birational 
geometry has been developed. For example see~\cite{Br1}, \cite{Ch}, 
\cite{Ka1}. 
Two smooth projective varieties $X$, $Y$ are called $K$-equivalent if 
and only if there is a common resolution $p\colon Z\to X$, $q\colon Z\to Y$
such that $p^{\ast}K_X =q^{\ast}K_Y$. 
Kawamata~\cite{Ka1} conjectured that 
derived categories are equivalent under $K$-equivalence. 
On the other hand Wang~\cite{Wan} conjectured that the deformation 
theories of complex structures are invariant under $K$-equivalence. 
Since derived category contains much information, 
it is reasonable to guess that 
 Kawamata's conjecture is stronger than Wang's conjecture. We 
will see the relationship between two conjectures using Theorem~\ref{mt}.
Recall that $X\stackrel{f}{\to}W\stackrel{g}{\gets}Y$ is called a 
flop if
\begin{itemize}
\item $f$ and $g$ are isomorphisms in codimension one. 
\item Relative Picard numbers of $f$, $g$ are one. 
\item $K_X =f^{\ast}K_W$, $K_Y =g^{\ast}K_W$. 
\item Birational map $g^{-1}\circ f \colon X\dashrightarrow Y$ is not an  
isomorphism.  
\end{itemize}
If $X$ and $Y$ are connected by flops, then $X$ and $Y$ are $K$-equivalent. 
We denote by $\kura (X)$ the Kuranishi deformation spaces, and by 
$T_0 \kura (X)$
its tangent space at the origin . 
Let $\xX \to \kura (X)$, $\yY \to \kura (Y)$ be 
Kuranishi families. For $\beta \in T_0 \kura (X)$, 
let $\xX _{\beta}$ be a scheme over 
$\mathbb{C}[\varepsilon]/(\varepsilon ^2)$, infinitesimal deformation of $X$
corresponding to $\beta$. 
\begin{thm}
Let $X$ and $Y$ be smooth projective varieties, which are 
connected by a flop, $X\stackrel{f}{\to}W\stackrel{g}{\gets}Y$.
Assume that there exists an object $\pP \in D^b(X\times Y)$, 
which is supported on $X\times _W Y$, such that the functor 
$\Phi _{X\to Y}^{\pP}\colon D^b(X)\to D^b(Y)$ gives an equivalence.
Then there exists an isomorphism $\phi _D \colon T_0 \kura (X)\to 
T_0 \kura (Y)$ such that $\Phi$ extends to an equivalence,
$$\Phi ^{\dag}\colon D^b(\Coh(\xX _{\beta}))
\longrightarrow D^b(\Coh(\yY _{\phi _D (\beta)})).$$
\end{thm}
\textit{Proof}.
Let $\phi _T \colon HT^2 (X)\to HT^2 (Y)$ be the isomorphism induced by 
$\Phi$. It suffices to show $\phi _T$ takes $(0,\beta,0)$ to $(0,\beta ',0)$.
Let $U\subset W$ be the maximum open subset on which $f$ and $g$ are 
isomorphic. Then, since $\codim (X\setminus U)\ge 2$, $\codim (Y\setminus U)
\ge 2$, and $f|_{X\setminus U}$, $g|_{Y\setminus U}$ has positive dimensional 
fibers, it follows that $\codim (W\setminus U)\ge 3$. On the other hand, since 
$\pP$ is supported on $X\times _W Y$, the following diagram commutes:
$$\begin{CD}
HT^2 (X) @>{\phi _T}>> HT^2 (Y) \\
@VVV  @VVV \\
HT^2 (U) @= HT^2 (U).
\end{CD}$$
Here the vertical arrows are restrictions.
Let $(\alpha ', \beta ', \gamma ')\cneq \phi _T (0, \beta, 0)$. 
By the above diagram, we have $\alpha ' |_{U}=0$, $\gamma '|_{U}=0$. 
It is clear $\gamma '=0$. On the other hand, 
since $\dR g_{\ast}\oO _Y =\oO _W$, 
we have $H^2 (Y, \oO _Y)\cong H^2 (W, \oO _W)$. Since 
$\codim (W\setminus U)\ge 3$, the restriction $H^2 (W, \oO _W)\to 
H^2 (U,\oO _U)$ is injective by~\cite{Sch}. Therefore $\alpha '=0$. 
$\quad$ q.e.d

\end{document}